\documentclass[leqno,11pt]{article}
\usepackage{latexsym}
\usepackage{amssymb}
\usepackage{amsfonts}
\usepackage{amsmath}
\usepackage{color}
\usepackage{epsfig}

\makeatletter
\long\def\unmarkedfootnote#1{{\long\def\@makefntext##1{##1}\footnotetext{#1}}}
\makeatother

\newtheorem{definition}{Definition}[section]

\newtheorem{lemma}[definition]{Lemma}

\newtheorem{theorem}[definition]{Theorem}

\newtheorem{remark}[definition]{Remark}

\def\m2{|\Omega | /2}
\def\M2{\frac{|\Omega |}{2}}
\def\u+{u_+^*}

\def\-p{\overline{p}}

\def\w0{{W_0^{1,p}(\Omega)}}

\def\R{\mathbb R}
\def\N{\mathbb N}

\def\rn{{{\R}^n}}

\newcommand{\hh}{{\cal H}^{n-1}}
\newcommand{\hu}{{\cal H}^{1}}
\newcommand{\medint}{-\kern  -,395cm\int}
\newcommand{\medintinrigo}{-\kern  -,315cm\int}
\newcommand{\medelle}{-\kern  -,235cm L}
\newcommand{\medellenrigo}{-\kern  -,180cm L}
\newcommand{\qed}{\thinspace\null\nobreak\hfill
\hbox{\vbox{\kern-.2pt\hrule height.2pt
depth.2pt\kern-.2pt\kern-.2pt \hbox to1.8mm {\kern-.2pt\vrule
width.4pt \kern-.2pt\raise1.8mm\vbox to.2pt{} \lower0pt\vtop
to.2pt{}\hfil\kern-.2pt \vrule
width.4pt\kern-.2pt}\kern-.2pt\kern-.2pt \hrule height.2pt
depth.2pt \kern-.2pt}}\par\medbreak}

\title{Balls minimize  trace constants in $BV$
}  \numberwithin{equation}{section}

\author{
  Andrea Cianchi\\
 {\it Dipartimento di Matematica e Informatica \lq\lq U. Dini", Universit\`a di Firenze}\\ {\it Piazza Ghiberti
27, 50122 Firenze, Italy}
\bigskip
\\
Vincenzo Ferone \\
{\it Dipartimento di Matematica e Applicazioni ``R. Caccioppoli"}
\\{\it Universit\`a di Napoli ``Federico II"}
\\ {\it Complesso Monte S.Angelo, Via Cintia, 80126 Napoli, Italy}
\bigskip
\\
Carlo Nitsch \\
{\it Dipartimento di Matematica e Applicazioni ``R. Caccioppoli"}
\\{\it Universit\`a di Napoli ``Federico II"}
\\ {\it Complesso Monte S.Angelo, Via Cintia, 80126 Napoli, Italy}
 \bigskip
\\
Cristina Trombetti \\
{\it Dipartimento di Matematica e Applicazioni ``R. Caccioppoli"}
\\{\it Universit\`a di Napoli ``Federico II"}
\\ {\it Complesso Monte S.Angelo, Via Cintia, 80126 Napoli, Italy}
}

\date{}

\begin{document}
\maketitle
\begin{abstract}\noindent
Balls are shown to have the smallest optimal constant, among all
admissible Euclidean domains,  in Poincar\'e type boundary trace
inequalities  for functions of bounded variation with vanishing
median or mean value.
\end{abstract}

\unmarkedfootnote {
}

\section{Introduction and main results}\label{sec1}

%
%
%
%We are concerned with boundary trace  inequalities, of Poincar\'e
%type, for
%%
%%
%% Consider the space $BV(\Omega )$ of
% functions of
%bounded variation in  open subsets $\Omega $ of $\rn$, $n \geq 2$.

A  branch of mathematical research
%at the edge between
which bridges  analysis and geometry is concerned with variational
problems for quantities of geometric--analytic nature associated
with  sets from some prescribed collection. Typically, the relevant
quantities are, in turn, expressed as a supremum or infimum of some
functional, defined on each set, which has often a physical meaning.
%
%
%is often associated with a physical quantity.
A prototypal  result in this area is the standard isoperimetric
inequality in the Euclidean space $\rn$. Further classical issues
amount to so called isoperimetric problems of mathematical physics,
and include inequalities for eigenvalues of elliptic operators and
isocapacitary inequalities. Most of these problems were originally
stated as conjectures. Some of them have been solved in the last
century via methods of the modern calculus of variations. Their
solution has led to such results as  Pol\'ya's theorem on De Saint
Venant's conjecture on cylindrical beams with the highest torsional
rigidity, Szeg\"o's theorem on Poincar\'e's conjecture on the body
of largest electrostatic capacity, Faber and Krahn's theorem on Lord
Rayleigh's conjecture on the lowest principal frequency of vibrating
clamped membranes. Other conjectures,
 including the
minimizing property of the ball for the first eigenvalue in a
fourth-order eigenvalue problem modeling the vibration of an elastic
 clamped plate (Szeg\"o conjecture), or the minimization property of the disk for the capacity in the family  of convex
 sets in the three-dimensional space with prescribed surface area (P\'olya-Szeg\"o conjecture),
 are still open, or are only known is special cases.
 %at least in their full generality.
 We do not even attempt an exhaustive bibliography on these topics.
 Let us
 just refer to the monographs and surveys \cite{AB, GGS, H, Ka, K, Ta} for an
 account of results and techniques in this field.
% classical and more recent deve
\par A class of quantities associated with open sets in $\rn$, whose
maximization  has traditionally attracted the attention of
specialists in functional and geometric analysis, is that of the
sharp constants in Sobolev-Poincar\'e type inequalities. Of course,
in many instances these constants can  be interpreted as eigenvalues
of an associated Euler equation. An overview of results and problems
in this connection can be found e.g. in \cite{BrV}. The
present contribution falls within this  line of investigations, and
focuses a minimization problem for the optimal constants in
Poincar\'e type inequalities for functions of bounded variation.
\par
Assume that $\Omega$ is a domain, namely a bounded connected open
set  in  $\rn$, $n \geq 2$. It is well known that if the boundary
$\partial \Omega$ of $\Omega$ is sufficiently regular, then a linear
operator if defined on the space $BV(\Omega )$ of functions of
bounded variation in $\Omega$, which associates with any function $u
\in BV(\Omega )$ its (suitably defined) boundary trace $\widetilde u
\in L^1(\partial \Omega )$. Here, $L^1(\partial \Omega )$ denotes
the Lebesgue space of integrable functions on $\partial \Omega$ with
respect to the $(n-1)$-dimensional Hausdorff measure $\hh$.
Moreover,  there exists a constant $C$, depending on $\Omega$, such
that
\begin{equation}\label{inf}
\inf _{c\in \R} \,\|\widetilde u - c\|_{L^1(\partial \Omega )}
%
%\int_{\partial\Omega} |\widetilde u - c| \, d {\mathcal
%H}^{n-1}(x)
%\|\widetilde u -{\rm med}_{\partial \Omega }\widetilde u
%\|_{L^1(\partial \Omega )}
 \le
 C \|Du\|(\Omega)
\end{equation}
for every $u \in BV(\Omega )$, where $\|Du\|(\Omega)$ stands for the
total variation  over $\Omega$ of the total variation of the
distributional gradient $Du$ of $u$ \cite[Theorem 9.6.4]{Mazbook}.
\par
A  property of $L^1$ norms ensures that the infimum in \eqref{inf}
is attained when $c$ agrees with a median of $\widetilde u$ on
$\partial \Omega$, given by
$${\rm med}_{\partial \Omega }\widetilde u = \sup \{t \in \mathbb
R: \hh (\{\widetilde u > t\}) > \hh (\partial \Omega )/2\}$$ (see
e.g. \cite[Lemma 3.1]{CianchiPick})
 Thus,
inequality \eqref{inf} is  equivalent to
\begin{equation}\label{trace}
%\inf _c \|\widetilde u - c\|_{L^1(\partial \Omega )}
%
%\int_{\partial\Omega} |\widetilde u - c| \, d {\mathcal
%H}^{n-1}(x)
\|\widetilde u -{\rm med}_{\partial \Omega }\widetilde u
\|_{L^1(\partial \Omega )}
 \le C_{\rm med}(\Omega) \|Du\|(\Omega)
\end{equation}
for every $u \in BV(\Omega )$, where  $C_{\rm med}(\Omega)$ denotes
the optimal -- smallest possible -- constant which renders
\eqref{trace} true.
\par
An other customary  Poincar\'e type trace inequality holds, when
${\rm med}_{\partial \Omega }\widetilde u$ is replaced with the mean
value $\widetilde u _{\partial \Omega }$ of $\widetilde u $ over
$\partial \Omega$, defined as
$$\widetilde u _{\partial \Omega } = \frac 1{\hh (\partial \Omega )}
\int _{\partial \Omega } \widetilde u \, d\hh (x)\,.$$
 The relevant inequality reads
\begin{equation}\label{tracemean}
%\int_{\partial\Omega} |\widetilde u - \widetilde u_{\partial
%\Omega }| \,d{\mathcal H}^{n-1}
\|\widetilde u - \widetilde u_{\partial \Omega }\|_{L^1(\partial
\Omega )} \le C_{\rm mv}(\Omega) \|Du \|(\Omega)
\end{equation}
for every $u \in BV(\Omega )$, where we have denoted by $C_{\rm
mv}(\Omega)$ the optimal constant in \eqref{tracemean}.
\par\noindent
Observe that, in the light of the above discussion, one has that
\begin{equation}\label{medmv}
C_{\rm med}(\Omega) \leq C_{\rm mv}(\Omega)
\end{equation}
for every domain $\Omega$. Also, note that both $C_{\rm
med}(\Omega)$ and $C_{\rm mv}(\Omega)$
%functionals
%in \eqref{minmax}
%\eqref{trace} and
%\eqref{tracemean}
are  invariant under dilations of $\Omega$,  and hence they only
depend on the shape of  $\Omega$, but not on its size.
\par
A minimal regularity assumption for inequalities \eqref{trace} and
\eqref{tracemean} to  hold is that  $\Omega $ be an admissible
domain, in the sense that
 $\hh (\partial \Omega
 )<\infty $,
 $\hh (\partial \Omega
 \setminus \partial ^M\Omega)=0$, and
 \begin{equation}\label{min}
%\sup _{E \subset \Omega} \frac{
\min\{\hh (\partial ^M E \cap
\partial \Omega) \, , \hh ( \partial \Omega \setminus \partial ^M E)\} \leq C \hh
(\partial ^M E \cap \Omega)
\end{equation}
for some positive constant $C$ and every measurable set $E \subset
\Omega$ \cite[Theorem 9.6.4]{Mazbook}. Here,  $\partial ^M $ denotes
the subset of the topological boundary, called the essential
boundary in  geometric measure theory. A local version of
\eqref{min} is, in fact, a necessary condition for the trace of $BV$
functions to be well defined on $\partial \Omega$ \cite{AG}.
%Note that, by \eqref{min}, any admissible domain is, in particular,
%connected.
Standard instances of admissible domains are provided by the
Lipschitz domains, namely bounded open sets whose boundary is
locally the graph of a Lipschitz function of $(n-1)$ variables.
\par
%We are concerned with the problem of
%
%The present paper addresses the problem of exhibiting the smallest
%possible constants in  boundary trace  inequalities, of Poincar\'e
%type, for
%%
%%
%% Consider the space $BV(\Omega )$ of
% functions of
%bounded variation in  open sets $\Omega \subset \rn$, $n \geq 2$, as
% $\Omega$ ranges among all admissible domains.
%
In this paper we address the problem of minimizing the trace
constants $C_{\rm med}(\Omega)$ and $C_{\rm mv}(\Omega)$, as
$\Omega$ ranges in the class of  all admissible domains $\Omega$ in
$\rn$.
 %Note
%that both $C_{\rm med}(\Omega)$ and $C_{\rm mv}(\Omega)$
%%functionals
%%in \eqref{minmax}
%%\eqref{trace} and
%%\eqref{tracemean}
%are  invariant under dilations of $\Omega$,  and hence the solutions
%to \eqref{minmax} account for the shape of the minimizer $\Omega$,
%but not for its size.
Heuristically
 speaking, domains with stretched shapes,
 %geometric configurations,
 such
 as sharp outward peaks or narrow passages, tend to have large
  values of $C_{\rm med}(\Omega)$ and $C_{\rm mv}(\Omega)$. One is
  thus  led to guess  that these constants attain their
  minimum value when $\Omega$ is a ball, in a sense the most rounded
    domain.
% Note that the  values of $C_{\rm med}(\Omega)$ and $C_{\rm
 % mv}(\Omega)$ when $\Omega$ is a ball are explicitly  exhibited in \cite[???]{Mazbook} and
  %\cite[???]{cianchitrace}, respectively.
%
%  least possible stretched
%  domain: a ball.
\par
In the two-dimensional case, the minimum problem for $C_{\rm
med}(\Omega)$   also arises in connection with
 questions of different nature. The   minimizing
property of the disk for $C_{\rm med}(\Omega)$  in classes of
admissible domains is known, and has been independently established
in \cite{KS, Es2, EGK}.
%A proof showing that the disk minimizes $C_{\rm
%med}(\Omega)$ in the class of planar convex domains is given in
%\cite{EGK}.
\par The higher-dimensional case appears to be open in the
existing literature.
%One difficulty in this kind of problems stems
%from the fact that they amount to $min$-$max$ problems, for
%functionals depending both on $\Omega$ and on $u$, of the form
%\begin{equation}\label{minmax}\inf _{\Omega} \sup _{u} \frac{\|\widetilde u -{\rm med}_{\partial \Omega
%}\widetilde u \|_{L^1(\partial \Omega )}}{ \|Du\|(\Omega)}\quad
%\hbox{and} \quad
% \inf _{\Omega} \sup _{u }
%\frac{\|\widetilde u -\widetilde u _{\partial \Omega }
%\|_{L^1(\partial \Omega )}}{ \|Du\|(\Omega)},
%\end{equation}
%respectively, as $\Omega$ ranges among all admissible domains in
%$\rn$, and $u$ among all non-constant functions in $BV(\Omega)$.
A lower estimate for $C_{\rm med}(\Omega)$, when $n \geq 3$, is
given in \cite[Theorem 6]{Es2}. This estimate, however, depends on
the geometry of $\partial \Omega$, and seems not to yield the
solution to the minimum problem for $C_{\rm med}(\Omega)$ in any
obvious way.
\par
Our  results confirm the above guess in any
dimension $n$, and also point out a singular phenomenon as far as
the uniqueness of minimizers is concerned.
%
%
%Our  results provide a solution to the problem at hand in any
%dimension $n$, and show that balls actually minimize both trace
%constants $C_{\rm med}(\Omega)$ and $C_{\rm mv}(\Omega)$. They also
%tell us that balls are the only minimizers, except for $C_{\rm
%mv}(\Omega)$ in the two-dimensional case, where a new non-uniqueness
%phenomenon occurs.
%
%
% point out an
%interesting difference between the situation for $C_{\rm
%med}(\Omega)$ and $C_{\rm mv}(\Omega)$ in the two-dimensional case.
\par
 Let us first consider $C_{\rm
med}(\Omega)$. In this regard, we have that the ball is the only
minimizer for $C_{\rm med}(\Omega)$  in any dimension $n \geq 2$.

\begin{theorem}\label{estimate1}
Let $\Omega$ be an admissible domain in $\rn$, $n \geq 2$.
% and let $B$ be a
%ball in $\rn$.
Then
%
%
%Under the hypotheses of Theorem \ref{mazya}
\begin{equation}\label{medianinequality}
C_{\rm med}(\Omega)\ge \sqrt{\pi} \,\frac n2 \frac{\Gamma (\frac
{n+1}2)} {\Gamma (\frac {n+2}2)}.
%K(B)
  %\frac{\omega_n}{2v_{n-1}}.
\end{equation}
%\textcolor{red}
Moreover, equality holds in \eqref{medianinequality}  if and only if
$\Omega$ is equivalent to a ball, up to a set of $\hh$ measure zero.
\end{theorem}

Our approach to Theorem \ref{estimate1} relies upon a
characterization of $C_{\rm med}(\Omega)$ as a genuinely geometric
quantity associated with $\Omega$, namely the optimal constant $C$
in \eqref{min}. Indeed, \cite[Theorem 9.5.2]{Mazbook} tells us that
\begin{equation}\label{admiss}
 C_{\rm med}(\Omega)=\sup_{E\subset\Omega}\frac{
\min\{\hh (\partial ^M E \cap
\partial \Omega) \, , \hh ( \partial \Omega \setminus \partial ^M
E)\}}{\hh (\partial ^M E \cap \Omega)},
%
% K(\Omega)=\sup_{E\subset\Omega}\frac{\min\{P_{C\Omega}
%(E), P_{C\Omega}(\Omega\setminus E)\}}{P_{\Omega} (E)}.
\end{equation}
where the supremum is extended over all measurable sets $E \subset
\Omega$ with positive Lebesgue measure.
\par\noindent
Note that, in particular, the constant appearing on the right-hand
side of \eqref{medianinequality} equals $\frac{n \omega _n}{2 \omega
_{n-1}}$, where $\omega _n=   \pi^{\frac n2}/\Gamma (1+\frac n2)$,
the Lebesgue measure of the unit ball in $\rn$. The fact that
 $C_{\rm
med}(\Omega)=\frac{n \omega _n}{2 \omega _{n-1}}$ when $\Omega $ is
a ball was proved in \cite{BM}, \cite[Theorem 9.5.2 and Corollary
9.4.4/3]{Mazbook} (see also \cite{BokSper} and \cite{Es2}). The
supremum in \eqref{admiss} is    attained at   a half-ball in this
case. Moreover, characteristic functions of half-balls  yield
equality in \eqref{trace}.
%
%
%
%In the case when $\Omega$ is a ball, the constant $C_{\rm
%med}(\Omega)$ is known to agree with $\frac{n \omega _n}{2 \omega
%_{n-1}}$, and  the supremum in \eqref{admiss} is    attained at   a
%half-ball \cite{BM}, \cite[Theorem 9.5.2 and Corollary
%9.4.4/3]{Mazbook} (see also \cite{BokSper} and \cite{Es2}).
%Functions which yield equality in \eqref{trace} are characteristic
%functions of half-balls.

\smallskip

\par
 We now take into account $C_{\rm mv}(\Omega)$. The next result shows
that the ball minimizes $C_{\rm mv}(\Omega)$  as well, and it is the
unique minimizer provided
  that $n \geq 3$. Interestingly enough, unlike   $C_{\rm
  med}(\Omega)$,   disks are not the only
minimizers of $C_{\rm mv}(\Omega)$ if $n=2$.

%
%\begin{theorem}\label{meanball}{\bf \cite[Theorem
%1.2]{cianchitrace}} Let $B$ be a ball in $\rn$. Then
%\begin{equation}\label{meanball1}
% K(B) =
%\begin{cases}  \frac{\omega_n}{2v_{n-1}} & \hbox{if $n \geq
%3$,}
%\\ 2  & \hbox{if $n =2$.}
%\end{cases}
%\end{equation}
%\end{theorem}

\begin{theorem}\label{estimatemean} Let $\Omega$ be an admissible
domain in $\rn$. If $n\ge 3$, then
\begin{equation}\label{estimatemean1}
%K(\Omega)\ge
%\begin{cases}  \frac{\omega_n}{2v_{n-1}} & \hbox{if $n \geq
%3$,}
%\\ 2  & \hbox{if $n =2$.}
%\end{cases}
C_{\rm mv}(\Omega)\ge \sqrt{\pi} \,\frac n2 \frac{\Gamma (\frac
{n+1}2)} {\Gamma (\frac {n+2}2)},
%
%\frac{\omega_n}{2v_{n-1}}
%\frac{\omega_n}{2v_{n-1}} \,\,(= K(B)).
\end{equation}
%\textcolor{blue}
and the equality holds in \eqref{estimatemean1}  if and  only if
$\Omega$ is equivalent to a ball, up to a set of $\hh$ measure zero.
\par\noindent
If $n=2$, then
\begin{equation}\label{estimatemean2}
C_{\rm mv}(\Omega)\ge 2,
\end{equation}
and the equality holds in \eqref{estimatemean2} if $\Omega$ is a
disk. However there exist open sets $\Omega$, that are not
equivalent to a disc,  for which equality yet holds in
\eqref{estimatemean2}.
\end{theorem}

Also the proof of Theorem \ref{estimatemean}  makes use of a
geometric characterization of $C_{\rm mv}(\Omega)$. This is provided
by \cite[Theorem 1.1]{cianchitrace}, and  reads
\begin{equation}\label{isopconstant}
C_{\rm mv}(\Omega) = \frac 2{\hh (\partial \Omega )} \sup _{E
\subset \Omega } \frac{\hh (\partial ^M E \cap
\partial \Omega ) \,\, \hh ( \partial \Omega \setminus \partial ^M E)}{\hh
(\partial ^M E \cap \Omega)},
\end{equation}
where the supremum is extended over all measurable sets $E \subset
\Omega$ with positive Lebesgue measure.
\par\noindent The value of $C_{\rm mv}(\Omega)$  when
$\Omega$ is a ball has recently been shown to coincide with the
right-hand side of either  \eqref{estimatemean1} or
\eqref{estimatemean2}, according to whether $n \geq 3$ or $n=2$
\cite[Theorem 1.2]{cianchitrace}. In the former case, $C_{\rm
mv}(\Omega)= C_{\rm med}(\Omega)$,  and the supremum in
\eqref{isopconstant} is achieved if $E$ is a half-ball. In the
latter case, however, $C_{\rm mv}(\Omega)> C_{\rm med}(\Omega)$, and
no maximizer exists on the right-hand side of \eqref{isopconstant}.
The supremum is approached along any sequence of circular segments
-- intersections of a disk with a half-plane -- whose measure tends
to zero. Accordingly, if $n \geq 3$, equality holds in
\eqref{tracemean} provided that $u$ is the characteristic function
of a half-ball, whereas, if $n=2$,  equality never holds in
\eqref{tracemean}, although the constant $C_{\rm mv}(\Omega)=2$ is
sharp, as demonstrated by any sequence of characteristic functions
of circular segments whose measure tends to zero. The lack of a
maximizer in \eqref{isopconstant} in the case when $\Omega$ is  a
disk allows, in a sense, for slight deformations of $\Omega$ which
do not affect $C_{\rm mv}(\Omega)$. As will be shown in the proof of
Theorem \ref{estimatemean}, a family  of domains $\Omega$ for which
$C_{\rm mv}(\Omega)$ agrees with
%
%is the same as
that of a disk consists in (nearly circular)
stadium-shaped sets.

\begin{remark}\label{rmed}
{\rm Inequalities \eqref{trace} and \eqref{tracemean} hold, in
particular, for every function $u$ in the Sobolev space
$W^{1,1}(\Omega) \subset BV(\Omega)$. Of course, in this case
$\|Du\|(\Omega)$   can be replaced with $\|\nabla u\|_{L^1(\Omega
)}$, where $\nabla u$ denotes the weak gradient of $u$.
%
%
%
%
%Boundary traces
%%on the boundary of  a Lipschitz domain $\Omega$
%%is,
%%$\partial \Omega$
%of functions $u$ from the Sobolev space $W^{1,1}(\Omega ) \subset
%BV(\Omega )$
%%, namely the Banach space of those functions $u \in
%%L^1(\Omega )$ whose first-order distributional gradient $\nabla u$
%%fulfills $|\nabla u| \in L^1(\Omega )$,
%are more classically defined on a Lipschits domain $\Omega$ as the
%limit of the standard traces on $\partial \Omega$ of approximating
%sequences
%%in $W^{1,1}(\Omega)$
%of smooth functions on $\overline \Omega$.
Let us emphasize that the constants $C_{\rm med}(\Omega)$ and
$C_{\rm mv}(\Omega)$ are   optimal  in the resulting Poincar\'e
trace inequalities in $W^{1,1}(\Omega)$ as well. Indeed,  any
function $u \in BV(\Omega)$ can be approximated by a sequence of
functions $u_k \in W^{1,1}(\Omega)$ in such a way that
$$\hbox{$\widetilde {u_k} = \widetilde u$ \quad and \quad $\lim _{k \to \infty}\|\nabla u_k\|_{L^1(\Omega )} = \|Du\|(\Omega
)$}.$$  The existence of the sequence $\{u_k\}$  follows, for
instance, from \cite[Theorem 1.17 and Remark 1.18]{Gi}.
%, or from
%\cite[Theorem 6]{AG}.
 Thus, Theorems \ref{estimate1} and \ref{estimatemean} also hold
if $C_{\rm med}(\Omega)$ and $C_{\rm mv}(\Omega)$  are interpreted
as the optimal constants in the trace inequalities \eqref{trace} and
\eqref{tracemean} for $u \in W^{1,1}(\Omega)$.}
\end{remark}

Optimal trace constants, and related shape optimization problems,
are the subject of various contributions, besides those already
mentioned above. Estimates for the constant $C$ in the Sobolev type
trace inequality
\begin{equation*}
\|\widetilde u \|_{L^1(\partial \Omega )}
 \le C (\|Du\|(\Omega) + \|u\|_{L^1(\Omega )})
\end{equation*}
for  $u \in BV(\Omega )$
 are
provided in \cite{AMR}.  The  inequality
\begin{equation}\label{mazya}\|u\|_{L^{\frac{n}{n-1}}(\Omega )} \leq \frac{\Gamma (1+ \frac
n2)^\frac 1n}{n \sqrt \pi} (\|Du\|(\Omega) +  \|\widetilde
u\|_{L^1(\partial \Omega )})
\end{equation} for $u \in W^{1,1}(\Omega)$, and hence for $u \in BV(\Omega )$,
where $\frac{\Gamma (1+ \frac n2)^\frac 1n}{n \sqrt \pi}$ is the
optimal constant, was proved in \cite{Matesi} (see also
\cite{Ma1969, Mazbook}). Versions of inequality \eqref{mazya} for
functions in the Sobolev space $W^{1,p}(\Omega)$, with $p>1$, can be
found in \cite{MV, MV1}. The paper \cite{Cimosertrace} contains a
 Poincar\'e trace inequality, with sharp exponential  constant,
for functions in the limiting Sobolev space $W^{1,n}(\Omega)$. The
optimal constant in the trace inequality for functions in
$W^{1,p}(\Omega)$, when $\Omega$ is a half-space,  was exhibited in
\cite{Es1} for
 $p=2$, and in \cite{Na} for any $p \in (1, n)$; the case $p=1$ is
 easy, as observed in \cite{cianchitrace}. A related Hardy-type trace
 inequality, with sharp constant, in a half-space is established in
 \cite{DDM}; an improved inequality, with remainder terms, is the
 object of \cite{AFV}. Related issues about optimal constants in Sobolev trace inequalities are
 discussed in \cite{BGP, Ro}.
 \par Let us finally mention that questions of a similar
 nature for mean-value Poincar\'e type inequalities  for functions in $BV(\Omega)$ and $W^{1,1}(\Omega)$, involving norms
 of $u$ in the whole of $\Omega$ instead of trace norms,  are
 treated in
 \cite{BoV, BrV, Cpoincare, EFKNT}. Contributions on optimal Poincar\'e inequalities in
 Sobolev spaces $W^{1,p}(\Omega)$, with $p>1$,   include  \cite{BK, DGS, DN, ENT, FNT, GW, Le}.

\section{A Cauchy formula for sets of finite perimeter}\label{sec2}

After recalling a few basic definitions and properties from
geometric measure theory, in this section we establish a version for
sets of finite perimeter of the classical Cauchy formula which
expresses the perimeter of an $n$-dimensional convex set in terms of
the measure of its $(n-1)$-dimensional projections.
\par
Let $E$ be a measurable set in $\rn$.  The upper and lower densities
$\overline D(E, x)$ and $\underline D(E, x)$ of E at a point $x \in
\rn$ are defined as $$ \overline D(E, x) = \limsup _{r\to 0} \frac
{\mathcal L^n(E \cap B_r(x))}{ \mathcal L^n(B_r(x))} \quad \hbox{
and} \quad \underline D(E, x) =\liminf _{r\to 0} \frac {\mathcal
L^n(E \cap B_r(x))}{ \mathcal L^n(B_r(x))},$$ respectively. Here,
$\mathcal L^n$ denotes the (outer) Lebesgue measure,
and $B_r(x)$ the ball centered at $x$, with radius $r$. When
$\overline D(E, x)$ and $\underline D(E, x)$ agree, their common
value is called the density of $E$ at $x$ and is denoted by $D(E,
x)$.  For each $\alpha \in [0, 1]$, the set $E^\alpha = \{x \in R^n
: D(E, x) = \alpha \}$ is called the set of points of density
$\alpha$ with respect to $E$, and is a Borel set. The set $E^1$ of
points of density $1$ with respect to $E$ agrees with $E$, up to
sets of Lebsegue measure zero. The essential boundary of $E$,
defined as
$$\partial ^M E = R^n \setminus (E^0 \cup E^1),$$ is also a Borel
set. Observe that $\partial ^M E \subset \partial E$.
\par\noindent
It is easily verified from the definition of essential boundary
that, if $E$ and $F$ are measurable subsets of $R^n$, then
\begin{equation}\label{boundaries}
\partial ^M (E \cup F) \cup \partial ^M (E \cap F) \subset
\partial ^M E \cup \partial ^M F.
\end{equation}
Note also that, if $E$ is any measurable set and $A$ is an open set,
then
\begin{equation}\label{boundaries1}
\partial ^M E \cap A \subset  \partial ^M (A \cap
E).
\end{equation}
Equation \eqref{boundaries1} follows from the  fact that, since $A$
is open,
$$\mathcal
L^n(E \cap B_r(x)) = \mathcal L^n((A \cap E )\cap B_r(x)) \quad
\hbox{if $x\in A$,}
$$
provided that $r$ is sufficiently small.
\par\noindent
Let $\Omega$ be an open set. If $E$ and $F$ are measurable subsets
of $\Omega$ such that $E\subset F$ (up to sets of zero Lebesgue
measure), then
\begin{equation}\label{boundaries0}
\partial ^M E \cap \partial ^M \Omega  \subset \partial ^M F \cap \partial ^M
\Omega.
\end{equation}
This is an easy consequence of the definition of essential boundary,
and of the fact that, if $x \in \partial ^M E \cap \partial ^M
\Omega$, then \begin{equation}  \limsup _{r\to 0} \frac {\mathcal
L^n(\Omega \cap B_r(x))}{ \mathcal L^n(B_r(x))} \geq \limsup _{r\to
0} \frac {\mathcal L^n(F \cap B_r(x))}{ \mathcal L^n(B_r(x))} \geq
\limsup _{r\to 0} \frac {\mathcal L^n(E \cap B_r(x))}{ \mathcal
L^n(B_r(x))}>0
\end{equation} and
\begin{equation}  1 >\liminf _{r\to 0} \frac {\mathcal
L^n(\Omega \cap B_r(x))}{ \mathcal L^n(B_r(x))} \geq \liminf _{r\to
0} \frac {\mathcal L^n(F \cap B_r(x))}{ \mathcal L^n(B_r(x))} \geq
\liminf _{r\to 0} \frac {\mathcal L^n(E \cap B_r(x))}{ \mathcal
L^n(B_r(x))}.
\end{equation}
As a consequence of \cite[Lemma 9.4.2]{Mazbook},
% $\hh (\partial ^M
%(\Omega \setminus E) \cap
%\partial ^M \Omega ) = \hh (\partial ^M
% \setminus \partial ^ME)$, or, equivalently,
\begin{equation}\label{boundaries4}
\hh (\partial ^M \Omega \cap \partial ^M (\Omega \setminus E)) = \hh
(\partial ^M \Omega
 \setminus \partial ^ME)
%
%\hh (\partial ^M E \cap \partial ^M \Omega ) + \hh (\partial ^M
%(\Omega \setminus E) \cap \partial ^M \Omega ) = P(\Omega )
\end{equation}
for every measurable subset $E$ of $\Omega$.
\par
The space $BV(\Omega)$ consists of  those  functions $u \in
L^1(\Omega)$ whose first-order distributional gradient $Du$ is a
vector-valued Radon measure with finite total variation
$\|Du\|(\Omega)$. The space $BV(\Omega )$ is a Banach space endowed
with the norm given by $\|u\|_{L^1(\Omega)} + \|Du\|(\Omega)$ for $u
\in BV(\Omega )$.
    \par\noindent
The boundary trace $\widetilde u $ of a function  $u \in BV(\Omega
)$
   can be
    %the function
%
%
%    and
%$u \in BV(\Omega)$, then the function
%$$\widetilde u :
%\partial \Omega \to \R\,$$
%is the trace of $u$ on $\partial \Omega$
defined for $\hh$-a.e. $x \in \partial \Omega$ as
\begin{equation}\label{trace0}
\widetilde u (x)= \lim _{r \to 0} \frac 1{\mathcal L ^n(B_r(x) \cap
\Omega)} \int _{B_r(x) \cap \Omega}u(y)\, dy\,,
\end{equation}
see \cite[Corollary 9.6.5]{Mazbook}. Note that this limit actually
 exists for $\hh$-a.e. $x \in
\partial \Omega$.  As recalled in Section \ref{sec1}, one has that $\widetilde u \in L^1(\partial \Omega
)$ for every function $u \in BV(\Omega)$. Moreover,
%
%for every  $u \in BV(\Omega)$, its trace $\widetilde u \in
%$L^1(\partial \Omega )$, the space of integrable functions on
%$\partial \Omega$ with respect to $\hh$,
 $L^1(\partial \Omega )$ cannot be replaced with any smaller Lebesgue
 space independent of $u$. Alternative definitions of
 the boundary trace of a function of bounded variation are available in the literature. One
  depends on the upper and lower approximate limits of the
extension of $u$ by $0$ outside $\Omega$ \cite[Definition
5.10.5]{Z}. Another one involves  the rough trace \cite[Section
9.5.1]{Mazbook}. Both of them
 coincide with $\widetilde u$, up
to subsets of $\partial \Omega$ of $\hh$-measure zero.
 \par\noindent
Traces
%on the boundary of  a Lipschitz domain $\Omega$
%is,
%$\partial \Omega$
of functions $u$ from the Sobolev space $W^{1,1}(\Omega )$
%, namely the Banach space of those functions $u \in
%L^1(\Omega )$ whose first-order distributional gradient $\nabla u$
%fulfills $|\nabla u| \in L^1(\Omega )$,
are more classically defined on the boundary of a Lipschitz domain
$\Omega$ as the limit of the restrictions to $\partial \Omega$ of
approximating sequences
%in $W^{1,1}(\Omega)$
of smooth functions on $\overline \Omega$.
%Since such an approximation is in
%general not possible for functions in $BV(\Omega)$, this
%definition of traces does not apply in $BV(\Omega)$, and one has
%to resort to \eqref{trace0}. However, for functions from
%$W^{1,1}(\Omega)$, the notions of traces in $W^{1,1}(\Omega)$ and
%in $BV(\Omega )$
This definition also yields a function on $\partial \Omega$ which
agrees with $\widetilde u$, up to subsets of $\partial \Omega$ of
$\hh$-measure zero.
\par
A measurable set $E\subset \rn$ is said to be of finite perimeter
 relative to  $\Omega$ if $D\chi _E $ is a
vector-valued Radon measure in $\Omega$ with finite total
variation in $\Omega$.
 The perimeter of $E$  relative to $\Omega$ is defined as
$$P(E; \Omega ) = \|D \chi _E
\|(\Omega ).$$ A result in geometric measure theory tells us that
$E$ is of finite perimeter in $\Omega$ if and only if $\hh (\partial
^M E \cap \Omega)< \infty$; moreover,
\begin{equation}\label{boundaries3}
P(E; \Omega ) = \hh (\partial ^M E \cap \Omega )
\end{equation}
%
%$\|D \chi _E \|(\Omega )= \hh (\partial ^M E \cap \Omega)$
\cite[Theorem 4.5.11]{Fed}. When $\Omega = \rn$, we denote $P(E;
\Omega )$ simply by $P(E)$, and call it the perimeter of $E$. Thus,
%$P(E)= \hh (\partial ^M E)$.
\begin{equation}\label{boundaries2}
P(E ) = \hh (\partial ^M E).
\end{equation}
Let $E$ be a set of finite perimeter in $\rn$. Then
the derivative $\nu ^E$ of the vector-valued measure $D \chi _E$
with respect to its total variation $|D \chi _E|$ exists, and
satifies $|\nu ^E (x)|=1$ for $\hh$-a.e. $x \in \partial ^ME$. The
vector $\nu ^E(x)$ is called the generalized inner normal to $E$ at
$x$.
%
%
%
%Let us denote by $P(\Omega )$ the perimeter of $\Omega$, and by
%$P(E; \Omega )$ the perimeter of a measurable set $E \subset \Omega$
%relative to $\Omega$. Then
%\begin{equation}\label{boundaries2}
%P(\Omega ) = \hh (\partial ^M \Omega ).
%\end{equation}
%Moreover,
%\begin{equation}\label{boundaries3}
%P(E; \Omega ) = \hh (\partial ^M E \cap \Omega )
%\end{equation}
%for every measurable subset $E$ of $\Omega$.
\par
Given $\nu \in \mathbb S^{n-1}$,  denote by $\nu ^\bot$ the
hyperplane which contains $0$ and is orthogonal to $\nu$. Given a
measurable set $E \subset \rn$, and $z \in \nu ^\bot$, we define
$$E_z^{\nu} =  \{r \in \R:  z+r\nu  \in E \}.$$
We also define the essential projection of $E$ on $\nu ^\bot$ as
$$\Pi _\nu (E)^+ = \{z \in \nu ^\bot : \mathcal L ^1 (E_z^{\nu})
>0\}.$$
Of course, the essential projection of $E$ agrees with its standard
projection if $E$ is open.
\par\noindent
 If,  $E$ is a bounded measurebale set, we set, for
$\nu \in \mathbb S^{n-1}$ and $z \in \Pi _\nu (E)^+$,
$$\phi _{E, \nu}(z) = \inf E_z^{\nu},$$
and, according to \cite[p. 233]{Gr}, we call the illuminated portion
of $E$ along $\nu$ the set
\begin{equation}\label{I}
I_\nu (E) = \{z+ \phi _{E, \nu}(z) \nu: \, z \in \Pi _\nu (E)^+\}.
\end{equation}

%
%$\Pi _\nu : \rn \to \nu ^\bot$ the orthogonal projection on $\nu
%^\bot$.

The  classical Cauchy formula tells us that, if $G$ is a convex set,
then
\begin{equation}\label{cauchy}
P(G)=\frac1{\omega _{n-1}}\int_{\mathbb S^{n-1}} {\mathcal
H}^{n-1}(\Pi_\nu (G )^+)\,d\nu\,,
\end{equation}
see \cite[Equation (32), Section 19]{BuZa} or \cite[Equation
(5.3.27)]{Schneider}. This formula was employed in the approach of
\cite{EGK}.
\par
A version of \eqref{cauchy} for  sets of finite perimeter is the
content of the following result.
%
%
%Our approach makes use of a version of this result for general sets
%of finite perimeter. This is the content of Theorem \ref{1} below.
%In what follows, $\omega _n$ denotes the Lebesgue measure of the
%unit ball in $\mathbb R^n$.

%Its statement and proof
%require a few notations and definitions from basic geometric measure
%theory.

\begin{theorem}\label{1}
Let $G$ be a set of finite perimeter and finite Lebesgue measure
in $\rn$.
%Denote by
%$\partial ^* \Omega $ the reduced boundary of $\Omega$.
 Then
\begin{equation}\label{cauchyper}
P(G)=\frac1{2\omega _{n-1}}\int_{\mathbb S^{n-1}} \bigg(\int _{\nu
^\bot} {\mathcal H}^{0}((\partial ^M G ) _z^\nu )\, d{\mathcal
H}^{n-1}(z )\,\bigg) d{\mathcal H}^{n-1}(\nu ).
\end{equation}
%where $\nu ^\bot$ denotes a hyperplane orthogonal to $\nu$.
In particular,
\begin{equation}\label{cauchyineq}
P(G)\ge\frac1{\omega _{n-1}}\int_{\mathbb S^{n-1}} {\mathcal
H}^{n-1}(\Pi_\nu (G )^+)\,d\nu\,.
\end{equation}
Moreover, the following facts are equivalent:
\par\noindent
(i) The equality holds in \eqref{cauchyineq};
\par\noindent
(ii) $G$ is equivalent to a convex set, up to sets of Lebesgue
measure zero;
\par\noindent
(iii) $G^1$ is convex.
%
%
%Equality holds in \eqref{cauchyineq} if and only if $G $ is
%equivalent to a convex  set.
\end{theorem}

\par\noindent
The discussion of the case of equality in \eqref{cauchyineq} in the
proof of Theorem \ref{1} makes use of the  next lemma, a slight
extension of a result of G. Alberti, reported  in the survey paper
\cite[Lemma 4.12]{fuscoisop}.

\begin{lemma}\label{dense}
Let $G$ be a measurable set in $\rn$, $n \geq 2$, such that, for
${\mathcal H}^{n-1}$- a.e. $\nu \in \mathbb S^{n-1}$, and for
${\mathcal H}^{n-1}$- a.e. $z \in \Pi_\nu (G )^+$,  the set $ G
_z^\nu$ is equivalent to an interval. Then $G^1$ is convex.
%
%
%
% having the property that
%\begin{equation}\label{dense0}
%\hbox{ for ${\mathcal H}^{n-1}$- a.e. $\nu \in \mathbb S^{n-1}$, and
%for ${\mathcal H}^{n-1}$- a.e. $z \in \Pi_\nu (G )^+$,  the set $ G
%_z^\nu$ equivalent to an interval.}
%\end{equation} Then the set of
%points of density $1$ with respect to $G$ is convex.
\end{lemma}
\par\noindent
{\bf Proof}.  Set $F=G^1$. Observe that, since $G$ and $F$ are
equivalent (up to sets of Lebesgue measure zero), then, as a
consequence of Fubini's theorem, the set $F$ satisfies the same
assumptions as $G$,
%
%property \eqref{dense0} is inherited by $F$,
namely
\begin{equation}\label{dense0F}
\hbox{ for ${\mathcal H}^{n-1}$- a.e. $\nu \in \mathbb S^{n-1}$, and
for ${\mathcal H}^{n-1}$- a.e. $z \in \Pi_\nu (F )^+$, the set $ F
_z^\nu$ is equivalent to an interval.} \end{equation}
\par\noindent Let $x_1, x_2 \in F$. Let $\nu \in \mathbb
S^{n-1}$, $\widehat z \in \nu ^\bot$, $y_1, y_2 \in \R$ be such that
$x_1 = \widehat z + y_1 \nu$, $x_2 = \widehat z + y_2 \nu$, with
$y_1 < y_2$. We have to show that any point $\widehat x$ of the form
$\widehat x = \widehat z + \widehat y \nu $, for some $\widehat y
\in (y_1, y_2)$, belongs to $F$.

Fix an orthogonal system whose $n$-th axis has the same direction
and orientation as $\nu$. Given $x \in \rn$ and $r>0$, denote by
$Q_r(x)$ the cube, centered at $x$ and with side-length $2r$, whose
sides are parallel to the coordinate axes of the relevant system.
Since $x_1, x_2$ have density $1$ with respect to $G$, they have
density $1$ with respect to $F$ as well. One has that $x \in F$ if
and only if $\frac{\mathcal L ^n (F\cap Q_r(x))}{2^n r^n} =1$ (see
e.g. \cite[Section 4.2]{fuscoisop}). Thus, for every $\varepsilon
>0$, there exists $r_\varepsilon >0$ such that, if $0 < r <
r_\varepsilon$, then
\begin{equation}    \label{dense1}
\frac{\mathcal L ^n (F\cap Q_r(x_i))}{2^n r^n} > 1- \varepsilon ,
\quad \qquad \hbox{for $i = 1,2$\,.}
\end{equation}
By property \eqref{dense0F},  there exists a sequence $\{\nu _k\}
\subset \mathbb S^{n-1}$, such that $\nu _k \to \nu$ as $k \to
\infty$, and   $ F _z^{\nu _k}$ is equivalent to an interval for
${\mathcal H}^{n-1}$-a.e. $z \in \Pi_{\nu _k} (F )^+$. From
\eqref{dense1} and  Fubini's theorem we deduce that
\begin{align} \label{dense2}
2^n r^n (1-\varepsilon ) & < \mathcal L ^n (F\cap Q_r(x_i))= \int
_{\Pi_{\nu _k}  (F \cap Q_r(x_i))^+} \mathcal L ^1 ((F \cap
Q_r(x_i))_z^{\nu _k}) \, d \mathcal H^{n-1}(z)  \\
\nonumber &  \leq 2r a_k \mathcal H^{n-1} (\Pi_{\nu _k}  (F \cap
Q_r(x_i))^+) \quad \qquad \hbox{for $i = 1,2$\,,}
\end{align}
for some sequence $\{a_k\}$ such that $a_k \geq 1$, and $a_k \to 1$
as $k \to \infty$. Hence,
\begin{equation}\label{dense3}
\mathcal H^{n-1} (\Pi_{\nu _k}  (F \cap Q_r(x_i))^+) \geq
\frac{2^{n-1} r^{n-1} (1-\varepsilon )}{a_k} \quad \hbox{for $i =
1,2$, and $k \in \mathbb N$\,.}
\end{equation}
Since $\Pi_{\nu _k}  (F \cap Q_r(x_1))^+)$ and $\Pi_{\nu _k}  (F
\cap Q_r(x_2))^+)$  satisfy \eqref{dense3}, and are contained in a
set which converges to an $(n-1)$-dimensional cube in $\nu ^\bot$ of
side-length $2r$ as $k \to \infty$, there exists another sequence
$\{b_k\}$ such that $b_k \geq 1$, and $b_k \to 1$ as $k \to \infty$,
such that
\begin{equation}\label{dense4}
\mathcal H^{n-1} (\Pi_{\nu _k}  (F \cap Q_r(x_1))^+ \cap \Pi_{\nu
_k}  (F \cap Q_r(x_2))^+) > \frac{2^{n-1} r^{n-1} (1-
2b_k\varepsilon )}{a_k}.
\end{equation}
Since $\{\nu_k\}$ is chosen in such a way that \eqref{dense0F} is
satisfied with $\nu = \nu_k$, we have that for ${\mathcal
H}^{n-1}$-a.e. $z \in \Pi_{\nu _k}  (F \cap Q_r(x_1))^+ \cap
\Pi_{\nu _k} (F \cap Q_r(x_2))^+$, the set $F_z^{\nu _k}$ is an
interval, and $\mathcal L^1((F \cap Q_r(x_i))_z^{\nu _k})>0$ for
$i=1, 2$. Thus, if $r$ is sufficiently small, depending on $y_1$,
$y_2$ and $\widehat y$, there exists a sequence
% sequences $\{c_k\}$ of
%numbers $c_k$ such that $c_k \leq 1$, $c_k \to 1$ as $k \to
%\infty$, and
$\{A_{k, r}\}$ of polyhedra in $\nu ^\bot _k$ such that
$$\mathcal H^{n-1} ((\Pi_{\nu _k}  (F \cap Q_r(x_1))^+ \cap
\Pi_{\nu _k}  (F \cap Q_r(x_2))^+) \setminus A_{k, r}) \to 0\quad
\hbox{ as $k \to \infty$,}$$ and
\begin{equation}\label{dense5}
\mathcal L^1 ((F \cap Q_r(\widehat x))_z^{\nu _k}) > 2r  \quad
\hbox{for $\mathcal H^{n-1}$ a.e. $z \in \Pi_{\nu _k}  (F \cap
Q_r(x_1))^+ \cap \Pi_{\nu _k}  (F \cap Q_r(x_2))^+ \cap A_{k, r}$.}
\end{equation}
Coupling \eqref{dense4} with \eqref{dense5} tells us that
\begin{equation}\label{dense6}
\mathcal L^n (F \cap  Q_r(\widehat x)) > 2r \Big(\frac{2^{n-1}
r^{n-1} (1- 2b_k\varepsilon )}{a_k} - c_{k,r}\Big)
\end{equation}
for some sequence $\{c_{k, r}\}$  such that $c_{k,r} \geq 0$ and
$c_{k, r} \to 0$ as $k \to \infty$ for every fixed (sufficiently
small) $r$. Passing to the limit in \eqref{dense6} as $k \to \infty$
yields
$$\frac{\mathcal L^n (F \cap  Q_r(\widehat x))}{2^n r^n} > (1- 2
\varepsilon),$$ whence
$$\liminf _{r \to 0}\frac{\mathcal L^n (F \cap  Q_r(\widehat x))}{2^n r^n} > (1- 2
\varepsilon).$$ Owing to the arbitrariness of $\varepsilon$, the
last equation ensures that $\widehat x$ has density $1$ with respect
to $F$, and hence $x \in F$. \qed

\par\noindent
{\bf Proof  of Theorem  \ref{1}}.
%Set $e_n = (0, \dots , 0, 1)$.
A special case of the coarea formula on rectifiable sets tells us
that
%\begin{equation}\label{coarea}\int _{\partial ^* G }|\nu
%_y^G (x)|\,d{\mathcal H}^{n-1}(x) = \int _{\mathbb R ^{n-1}}
%{\mathcal H}^{0}(\partial ^* G   \cap \Pi _{e_n} ^{-1}(x'))\,
%dx', \end{equation}
% where we have set $x=(x',
%y)$, $x' \in \mathbb R ^{n-1}$, $y \in \R$, and $\nu ^G (x)$
%denotes the generalized inner normal to $G$ at x, which is
%well defined at each point of $\partial ^* G$, and $\nu
%_y^G (x)$ stands for its component along $y$.
% \par\noindent
% Analogously,
 for each $\nu \in \mathbb S^{n-1}$,
\begin{equation}\label{coarea}\int _{\partial ^M G }|\nu ^G (x)\cdot \nu |\,d{\mathcal H}^{n-1}(x) =
\int _{\nu ^\bot} {\mathcal H}^{0}((\partial ^M G )_z^\nu )\,
d{\mathcal H}^{n-1}(z )
\end{equation}
(see, for instance, \cite[Theorem F]{CCF}).
%where  $\nu ^G (x)$ denotes the generalized inner normal to $G$
%at $x$, which is well defined at ${\mathcal H}^{n-1}$ a.e.  point $x
%\in \partial ^M G$.

On integrating equation \eqref{coarea} with respect to $\nu$ over
$\mathbb S^{n-1}$ yields
\begin{equation}\label{cauchyper3}\int _{\mathbb S^{n-1}}\int _{\partial ^M G }|\nu ^G (x)\cdot \nu |
\,d{\mathcal H}^{n-1}(x) d{\mathcal H}^{n-1}(\nu ) = \int
_{\mathbb S^{n-1}} \bigg(\int _{\nu ^\bot} {\mathcal
H}^{0}((\partial ^M G )_z^\nu )\, d{\mathcal H}^{n-1}(z ) \bigg)
d{\mathcal H}^{n-1}(\nu ).
\end{equation}
By Fubini's theorem
\begin{align}\label{cauchyper4}
\int _{\mathbb S^{n-1}}\int _{\partial ^M G }|\nu ^G (x)\cdot \nu
|\,d{\mathcal H}^{n-1}(x) d{\mathcal H}^{n-1}(\nu ) & =
\int_{\partial ^M G } \int _{\mathbb S^{n-1}}|\nu ^G (x)\cdot
\nu |\, d{\mathcal H}^{n-1}(\nu )\,d{\mathcal H}^{n-1}(x)\\
\nonumber &  = 2 \omega_{n-1} \int_{\partial ^M G } \,d{\mathcal
H}^{n-1}(x)
 = 2 \omega _{n-1} P(G ), \end{align} where the second
inequality holds since the integral
$$
\int _{\mathbb S^{n-1}}|\nu ^G (x)\cdot \nu |\, d{\mathcal
H}^{n-1}(\nu )$$ is clearly independent of $G$ and $x$, and equals
$2 \omega _{n-1}$. Equation \eqref{cauchyper} thus follows from
\eqref{cauchyper3} and \eqref{cauchyper4}.
\par\noindent
Let us now focus on \eqref{cauchyineq}. We claim that, for every
$\nu \in \mathbb S^{n-1}$,
 %Since $G$ is a set of finite perimeter (and $G
%\neq \rn$), one has that
\begin{equation}\label{dis2}
{\mathcal H}^{0}((\partial ^M G )_z^\nu ) \geq 2 \quad \hbox{for
$\mathcal H ^{n-1}$-a.e. $z \in \Pi_\nu (G)^+$.}
\end{equation}
Indeed,  for every $\nu \in \mathbb S^{n-1}$, there exists a Borel
subset $B_{G, \nu}$ of $\Pi_\nu (G )^+$ such that $\mathcal H ^{n-1}
(\Pi_\nu (G )^+ \setminus B_{G, \nu})=0$, and, for every $z \in
B_G$, the set $G _z^\nu $
%$$G _z = \{y \in \R:  z+ y\nu  \in G \}$$
is  of finite perimeter and measure in $\R$, and
 \begin{equation}\label{dic20}(\partial ^M G )_z^\nu = \partial ^M (G _z ^\nu)
 \end{equation}
 (see e.g.  \cite[Theorem G]{CCF}).
 Note that the fact that $\mathcal L^1 (G
_z^\nu) <\infty$ for $z \in B_{G, \nu}$ is a consequence of the
assumption  $\mathcal L ^n (G) < \infty$, since, by Fubini's
theorem,
$$\mathcal L ^n (G) = \int _{\Pi_\nu (G )^+} \mathcal L^1 (G
_z^\nu)\, dz.$$
 By the isoperimetric inequality in $\R$,
\begin{equation}\label{isopunid}{\mathcal H}^{0} (\partial ^M (G _z ^\nu)) \geq 2
 \end{equation}
for every $z \in B_{G, \nu}$, and hence for $\hh$- a.e. $z \in
\Pi_\nu (G )^+$. Moreover, the equality holds in \eqref{isopunid} if
and only if $G _z ^\nu$ is equivalent to an interval.
%Since, obviously,
%$${\mathcal H}^{0}(\partial ^MG \cap \Pi _\nu ^{-1} (z) ) = {\mathcal H}^{0}((\partial ^M G
%)_z) \quad \hbox{for $z \in \nu ^\bot$,}$$
Thus, \eqref{dis2} follows from \eqref{dic20} and \eqref{isopunid},
and one has that the equality holds in \eqref{dis2} if and only if $
G _z^\nu$ is equivalent to an interval.
\par\noindent
By \eqref{dis2}
\begin{equation}\label{cauchyper5}
\int _{\nu ^\bot} {\mathcal H}^{0}((\partial ^M G )_z^\nu )\,
d{\mathcal H}^{n-1}(z ) \geq 2 {\mathcal H}^{n-1}(\Pi_\nu
(G)^+)\qquad  \hbox {for every $\nu \in \mathbb S^{n-1}$.}
\end{equation}
Inequality \eqref{cauchyineq} is a consequence of
\eqref{cauchyper} and \eqref{cauchyper5}.
\par\noindent
As far as the case of equality in \eqref{cauchyineq} is concerned,
if (ii) holds, then the equality holds in inequalities
\eqref{dis2}--\eqref{cauchyper5}, and hence also in
\eqref{cauchyineq}, whence (i) follows. The fact that (iii) implies
(ii) is a consequence of the equivalence of $G$ and $G^1$ up to sets
of Lebesgue measure zero. It remains to show that (i) implies (iii).
Assume that (i) holds. Then
%
%
%
%Assume now that equality holds in \eqref{cauchyineq}. Then
equality holds in \eqref{cauchyper5} for ${\mathcal H}^{n-1}$- a.e.
$\nu \in \mathbb S^{n-1}$. Hence,  for ${\mathcal H}^{n-1}$- a.e.
$\nu \in \mathbb S^{n-1}$, equality also holds in \eqref{isopunid}
for ${\mathcal H}^{n-1}$ a.e.  $z \in \Pi_\nu (G )^+$. Thus, for
${\mathcal H}^{n-1}$- a.e. $\nu \in \mathbb S^{n-1}$, the set  $ G
_z^\nu$ equivalent to an interval   for ${\mathcal H}^{n-1}$- a.e.
$z \in \Pi_\nu (G )^+$. Property (iii) hence follows via Lemma
\ref{dense}.
%
%As a consequence of Lemma \ref{dense} below, the set of points of
%density $1$ with respect to $G$ is in fact convex, and hence $G$
%itself is equivalent to a convex set.
\qed

\section{Proofs of the main results}\label{sec3}

We begin by accomplishing the proof of Theorem \ref{estimate1}.
\medskip
\par\noindent
 {\bf Proof of Theorem \ref{estimate1}}.
%More generally, when $\Omega$ is a set of finite perimeter,
%then \eqref{cauchy} has to be replaced with the inequality
%\begin{equation}\label{cauchyineq}
%$$P(\Omega)\ge\frac1{v_{n-1}}\int_{\mathbb S^{n-1}} {\mathcal H}^{n-1}(\Pi_\nu\Omega)\,d\nu,$$
%\end{equation}
%(see Lemma \ref{1} below). Hence, the conclusion follows as above.
%then $H_\nu$ introduced above may not be defined for all $\nu\in\mathcal{S}^{n-1}$.
For any given $\nu\in\mathbb{S}^{n-1}$, there exist two sequences of
open half-spaces $\{H_{\nu ,i}^+\}_{i\in\N}$ and $\{H_{\nu
,i}^-\}_{i\in\N}$ such that, for $i \in \mathbb N$,
%
%
%on defining $H_\nu ^+=\displaystyle\cup_i H_{\nu ,i}^+$ and $H_\nu
%^-=\displaystyle\cup_i H_{\nu ,i}^-$:
%\begin{itemize}
%\item
\begin{align}\label{A}
& \hbox{$\pm\nu$ is the outer unit normal to $H_{\nu ,i}^\pm$ on
$\partial H_{\nu ,i}^\pm$,} \nonumber \\ \nonumber &  \hbox{$H_{\nu
,i}^\pm\subset H_{\nu , i+1}^\pm$,}\\ \nonumber % for $i \in \mathbb N$, &
& \hbox{$\overline{H_{\nu ,i}^+}\cap \overline{H_{\nu
,i}^-}=\emptyset$,}\\ & \hbox{ ${\mathcal
H}^{n-1}(\partial\Omega\cap\partial^M(H_{\nu ,i}^{\pm}\cap
\Omega))\le P(\Omega)/2$,}
\end{align}
and, on defining $H_\nu ^+=\displaystyle\cup_i H_{\nu ,i}^+$ and
$H_\nu ^-=\displaystyle\cup_i H_{\nu ,i}^-$,
$$\partial H_\nu ^+ =
\partial H_\nu^-.$$
 Notice that here, and in
similar occurrences below, the use of just $\partial \Omega$ instead
of $\partial ^M \Omega$ is allowed by the fact that $\Omega$ is an
admissible domain.
\par\noindent
 One has that
\begin{multline}\label{dic1}
{\mathcal H}^{n-1}(\partial\Omega\cap\partial^M(H_\nu^{+}\cap
\Omega)) + {\mathcal
H}^{n-1}(\partial\Omega\cap\partial^M(H_\nu^{-}\cap \Omega))
\\ {= {\mathcal H}^{n-1}(\partial\Omega\cap\partial^M(H_\nu^{+}\cap
\Omega)) + {\mathcal
H}^{n-1}(\partial\Omega\setminus\partial^M(H_\nu^{+}\cap \Omega))}=
P(\Omega).
\end{multline}
 Observe that the first equality in \eqref{dic1}  follows from
 \eqref{boundaries4} applied with $E= H_\nu^{+}\cap \Omega$, since
 $\Omega \setminus E = (H_\nu^{-}\cap \Omega) \cup (\partial H_\nu^{-}\cap
 \Omega)$, and hence  $\Omega \setminus E = (H_\nu^{-}\cap \Omega)$, up
 to sets of Lebesgue measure zero. The second equality in
 \eqref{dic1} holds since $\hh$ is a measure on Borel sets.
 \par\noindent
 {
 Let us define
 $$H_\nu = H_\nu ^+.$$
 We claim that
 \begin{equation}\label{B}
{\mathcal H}^{n-1}(\partial\Omega\cap\partial^M(H_\nu \cap
\Omega))=P(\Omega)/2 \quad \hbox{for $\hh$-a.e. $\nu \in \mathbb S
^{n-1}$.}
\end{equation}
Owing to \eqref{dic1}, equation \eqref{B} only fails if
%
%
%Now,  two cases may occur. If ${\mathcal
%H}^{n-1}(\partial\Omega\cap\partial^M(H^{\pm}\cap
%\Omega))=P(\Omega)/2$, then we define $H_\nu=H^+$ and set
%\begin{equation}\label{h} h(\nu)={\mathcal H}^{n-1}(\Omega\cap\partial^M(H_\nu\cap
%\Omega)).
%\end{equation}
%Consider next the case when
\begin{equation}\label{C}
{\mathcal H}^{n-1}(\partial\Omega\cap\partial^M(H_\nu^{+}\cap
\Omega)) \neq {\mathcal H}^{n-1}(\partial\Omega\cap\partial^M(H_\nu
^{-}\cap \Omega)).
\end{equation} Thus, in order to prove our claim, it suffices to
show that \eqref{C} can only hold for $\nu$ in a countable subset of
$\mathbb S ^{n-1}$. Assume that \eqref{C} is on force for some $\nu
\in \mathbb S ^{n-1}$.}
 Then
%$$\lim_i {\mathcal H}^{n-1}(\partial\Omega\cap \partial (H_{\nu ,i}^\pm\cap\Omega))\le P(\Omega)/2$$
\begin{eqnarray}\label{dic2}
&\lim_{i \to \infty} &\min\left\{ {\mathcal
H}^{n-1}(\partial\Omega\cap
\partial^M (H_{\nu ,i}^+\cap\Omega)),{\mathcal
H}^{n-1}(\partial\Omega\cap
\partial^M (H_{\nu ,i}^-\cap\Omega))\right\}\\ \nonumber &&\le \min\left\{
{\mathcal H}^{n-1}(\partial\Omega\cap \partial
^M(H_\nu^+\cap\Omega)),{\mathcal H}^{n-1}(\partial\Omega\cap
\partial ^M(H_\nu^-\cap\Omega))\right\}< P(\Omega)/2,
\end{eqnarray}
where  the first inequality  holds owing to equation
\eqref{boundaries0} with $E=H_{\nu ,i}^{\pm}\cap\Omega$ and
$F=H^{\pm}_\nu \cap \Omega$, {and the second one by
\eqref{A} and \eqref{C}}. The following chain holds:
 %By \eqref{dic2},
\begin{align}\label{dic3}
& P(\Omega)= {\mathcal H}^{n-1}(\partial\Omega\cap
H_\nu^+)+{\mathcal H}^{n-1}(\partial\Omega\cap H_\nu^-)+{\mathcal
H}^{n-1}(\partial\Omega\cap \partial H_\nu^+)\\ \nonumber &= \lim_{i
\to \infty}\left ({\mathcal H}^{n-1}(\partial\Omega\cap H_{\nu
,i}^+)+{\mathcal H}^{n-1}(\partial\Omega\cap H_{\nu
,i}^-)\right)+{\mathcal H}^{n-1}(\partial\Omega\cap \partial H_\nu^+)\\
\nonumber &\le \lim_{i \to \infty}\left ({\mathcal
H}^{n-1}(\partial\Omega\cap
\partial^M (H_{\nu ,i}^+\cap\Omega))+{\mathcal
H}^{n-1}(\partial\Omega\cap
\partial^M (H_{\nu ,i}^-\cap\Omega))\right)+{\mathcal
H}^{n-1}(\partial\Omega\cap
\partial H_\nu^+)\\ \nonumber
&\le   \lim_{i \to \infty}\min\left\{ {\mathcal
H}^{n-1}(\partial\Omega\cap
\partial^M (H_{\nu ,i}^+\cap\Omega),{\mathcal H}^{n-1}(\partial\Omega\cap
\partial^M (H_{\nu ,i}^-\cap\Omega)\right\} + P(\Omega)/2 + {\mathcal
H}^{n-1}(\partial\Omega\cap \partial H_\nu^+)\\ \nonumber &<
P(\Omega) + {\mathcal H}^{n-1}(\partial\Omega\cap \partial H_\nu^+).
\end{align}
Note that the first equality in \eqref{dic3} holds since $\hh$ is a
measure when restricted to Borel sets, and hence its is (countably)
additive on disjoint Borel sets, the second equality again relies
upon the fact that $\hh$ is a measure on Borel sets, and
$\partial\Omega\cap H_{\nu ,i}^{\pm} \nearrow \partial\Omega\cap
H_\nu ^{\pm}$, the first inequality is a consequence of the
inclusion
\begin{equation}\label{dic10}
\partial^M\Omega\cap H_{\nu ,i}^{\pm} \subset \partial^M\Omega\cap
\partial^M
(H_{\nu ,i}^{\pm }\cap\Omega),
\end{equation}
 which, in turn, follows from
\eqref{boundaries1} applied with $A=H_{\nu ,i}^{\pm}$,
{and the second and third inequality are consequences
of \eqref{A} and \eqref{dic2}, respectively.} Equation \eqref{dic3}
implies that ${\mathcal H}^{n-1}(\partial\Omega\cap
\partial H_\nu ^+)>0$. {Since $\Omega$ ha finite
perimeter, the latter inequality can  hold at most for $\nu$ in a
countable subset of $\mathbb {S}^{n-1}$.  Hence, our claim follows.}
%
%
%, and hence  ${\mathcal
%H}^{n-1}(\partial\Omega\cap\partial^M(H^{+}\cap \Omega)) \neq
%{\mathcal H}^{n-1}(\partial\Omega\cap\partial^M(H^{-}\cap \Omega))$
% at most for $\nu$ in a countable subset of $\mathbb {S}^{n-1}$. It thus turns
%out that the half space $H_\nu$ and the function $h(\nu)$ are
%well-defined as above for $\hh$-a.e. $\nu \in \mathbb {S}^{n-1}$.
\par\noindent
Now, by \eqref{boundaries}, $\partial ^M (H_\nu \cap \Omega )
\subset \partial ^M H_\nu \cup \partial ^M \Omega=
\partial  H_\nu \cup \partial ^M \Omega $. Thus, $\Omega \cap \partial ^M
(H_\nu \cap \Omega ) \subset \Omega \cap (\partial  H_\nu \cup
\partial ^M \Omega) = (\partial  H_\nu \cap \Omega ) \cup
(\partial ^M \Omega \cap \Omega ) =
\partial H_\nu \cap \Omega$. Hence, $$\hh(\Omega \cap \partial ^M
(H_\nu \cap \Omega ) ) \leq \hh (\partial H_\nu \cap \Omega).$$
Since $\nu$ is orthogonal to the hyperplane $\partial H_\nu$,
$$\hh (\partial H_\nu \cap \Omega) \leq \hh(\Pi _\nu (\Omega )^+).$$
Altogether, we obtain $${\mathcal H}^{n-1}(\Pi_\nu(\Omega)^+)\ge
{{\mathcal H}^{n-1}(\Omega\cap\partial^M(H_\nu\cap
\Omega))} \quad \hbox{for $\hh$-a.e. $\nu \in \mathbb {S}^{n-1}$.}$$
Hence, by \eqref{cauchyineq},
\begin{equation}\label{cauchy1}
P(\Omega)\ge \frac1{\omega _{n-1}}\int_{\mathbb S^{n-1}} {\mathcal
H}^{n-1}(\Pi_\nu(\Omega)^+)\,d\nu \geq \frac1{\omega
_{n-1}}\int_{\mathbb S^{n-1}} {{\mathcal
H}^{n-1}(\Omega\cap\partial^M(H_\nu\cap \Omega))}\,d\nu.
\end{equation}
Finally, from the definition of $C_{\rm med}(\Omega)$ and \eqref{cauchy1}
 %observing that \eqref{cauchy1} is still in force,
we have that
\begin{eqnarray}\label{C1}C_{\rm med}(\Omega)&\ge& \frac{P(\Omega)}2\frac1{\displaystyle\inf_{H \text{ half-space}\atop
{\mathcal H}^{n-1}(\partial\Omega\cap\partial^M(H\cap
\Omega))=P(\Omega)/2} {\mathcal H}^{n-1}(\Omega\cap\partial^M(H\cap
\Omega))}\\ \nonumber &\ge&
\frac{P(\Omega)}2\frac1{\displaystyle\frac1{n\omega_{n}}\int_{\mathbb
S^{n-1}} {{\mathcal
H}^{n-1}(\Omega\cap\partial^M(H_\nu\cap
\Omega))}\,d\nu}\ge\frac{n\omega_n}{2\omega _{n-1}} = \sqrt{\pi}
\,\frac n2 \frac{\Gamma (\frac {n+1}2)} {\Gamma (\frac {n+2}2)},
\end{eqnarray}
whence \eqref{medianinequality} follows.
\par
As for the equality case in \eqref{medianinequality}, observe that
if the equality  holds in \eqref{medianinequality}, then it also
holds in the chain of inequalities \eqref{C1}, and hence  in
\eqref{cauchy1} as well. Hence, by Theorem \ref{1},
 $\Omega ^1$
 % (the set of points of density 1 with respect to
%$\Omega$)
 is a convex set. The convexity of $\Omega ^1$ implies that
it is an open set. To verify this assertion, it suffices to show
that  $\partial\Omega ^1\cap \Omega ^1=\emptyset$. The latter
equality is in turn a consequence of the fact that, since $\Omega
^1$ is convex, if $x\in\partial\Omega ^1$ then
$$D(\Omega ^1, x)=\overline D(\Omega ^1, x)=\underline D(\Omega ^1, x)\in(0,1),$$
and hence $\partial \Omega ^1 = \partial ^M \Omega ^1 = \partial ^M
\Omega$. The openness of $\Omega$ implies that
%
%
%
% $\mathcal L^n(\Omega \cap \Omega ^1)=0$
%$$ D(\Omega, x)=D(\Omega ^1, x).$$
% \textcolor{blue}{  Since $\Omega$ is open, we have that
$\Omega\subset\Omega
 ^1$. We claim that
\begin{equation}\label{dic8}
 {\mathcal H}^{n-1}(\Omega
 ^1\setminus \Omega)=0.
 \end{equation}
  %The fact that $\Omega$ is admissible
% entails that
% \begin{equation}\label{dic4}
% \partial \Omega = \partial ^M \Omega \quad \hbox{up to sets of
% $\hh$ measure zero}.
% \end{equation}
% Owing to the convexity of $\Omega ^1$, we also have that
%\begin{equation}\label{dic5}
% \partial \Omega^1 = \partial ^M \Omega ^1.
% \end{equation}
%On the other hand, since $\Omega = \Omega ^1$ up to sets of Lebesgue
%measure zero,
%\begin{equation}\label{dic6}
% \partial ^M \Omega  = \partial ^M \Omega ^1.
% \end{equation}
%From \eqref{dic4}--\eqref{dic6} we obtain that
%\begin{equation}\label{dic7}
% \hh (\partial \Omega \bigtriangleup \partial  \Omega ^1 ) =0.
% \end{equation}
Indeed, since $\Omega \subset \Omega ^1 \subset \overline {\Omega
}$, we have that $\overline {\Omega ^1} = \overline {\Omega }$.
Thus, $\Omega
 ^1\setminus \Omega \subset \partial \Omega $. Inasmuch as  $\Omega
 ^1 \cap \partial ^M \Omega = \emptyset$, one has $\Omega ^1\setminus \Omega  \subset
 \partial \Omega \setminus \partial ^M \Omega$,  and hence \eqref{dic8}
 follows since $\Omega$ is admissible.
 \par\noindent
 It only remains to show that $\Omega ^1$ is a ball. Let
 $H_\nu$ denote the half-space defined as above, with $\Omega$
 replaced with $\Omega ^1$. The convexity of $\Omega ^1$ ensures
 that, for every $\nu\in\mathbb S^{n-1}$, $H_\nu$ agrees with the
unique open half space such that $\nu $ is the outer normal  to
$\partial H_\nu$, and
\begin{equation}\label{Hnu}
{\hh (\partial \Omega ^1 \cap H_\nu)} = {\mathcal
H}^{n-1}(\partial \Omega ^1 \cap
\partial^M(H_\nu\cap \Omega ^1))=P(\Omega ^1)/2.\end{equation}
 Since the
equality holds in the chain of inequalities \eqref{cauchy1} and
\eqref{C1},
\begin{equation}\label{equality6}
{\mathcal H}^{n-1}(\Omega ^1\cap\partial^M(H_\nu\cap \Omega
^1))={\mathcal H}^{n-1}(\Pi_\nu(\Omega
^1)^+)=P(\Omega)\frac{\omega_{n-1}}{n\omega_{n}},\qquad\hbox{for every
}\nu\in\mathbb S^{n-1}.
\end{equation}
Note that \eqref{equality6} actually holds for every $\nu\in\mathbb
S^{n-1}$, since ${\mathcal H}^{n-1}(\Omega
^1\cap\partial^M(H_\nu\cap \Omega ^1)$ and ${\mathcal
H}^{n-1}(\Pi_\nu(\Omega ^1)^+)$ are continuous functions of $\nu$,
owing to the convexity of $\Omega ^1$.
\par\noindent
Our next step consists in proving that $\Omega ^1$ is, in fact,
strictly convex.  To this purpose, it suffices to show that no
segment is contained in $\partial \Omega ^1$. Let us
 assume, by contradiction,  that there
exist a straight line intersecting $\partial \Omega^1$ in a whole
segment $\Sigma$. Denote by
%$x_1, x_2\in\partial\Omega^1 $ such
%that the whole segment $\Sigma$ with endpoints
$x_1$ and $x_2$ the endpoints of $\Sigma$.
%
% contained in $\partial \Omega^1$.
%
%$$
%\textcolor{red} {t\,x_1+(1-t)\,x_2 \in \partial \Omega^1 \quad\mbox{for all $t \in[0,1]$} }
%$$
%
%\textcolor{red}{
Set $\nu=\frac{x_1-x_2}{|x_1-x_2|}$.
 %and by
%$$s=\left\{ t\, x_1+(1-t)\, x_2: \; t\in[0,1]\right\}.$$
Observe that $\partial H_\nu\cap \Sigma\ne\emptyset$, otherwise $\hh
(\Pi_\nu(\Omega^1)^+)>\hh(\Omega^1\cap\partial^M(H_\nu\cap\Omega^1))$,
thus contradicting \eqref{equality6}. Let  $\widehat y$ be the point
such that   $\partial H_\nu\cap \Sigma=\{\widehat y\}$, and let
$\widehat z$ be any point in  $\Sigma$, different from $x_1$, $x_2$,
and $ \widehat y$. Let $H$ be an open half-space such that $\Sigma
\subset
\partial H$ and $\Omega ^1 \subset H$ (in particular, $ \partial H$
is a support hyperplane to the convex set $\Omega ^1$), and let $\mu
\in \mathbb S^{n-1}$ be the outward unit normal vector to $H$ on
$\partial H$. One has that
%
%
%Define $\mathcal{C}=\text{span}(\mu,\nu) \cap \mathbb{S}^{n-1}$, and
%note that
$${\text{span}(\mu,\nu) \cap \mathbb{S}^{n-1} = \{\xi (\vartheta ): \xi (\vartheta )= \cos (\vartheta ) \mu + \sin
(\vartheta ) \nu \,\, \hbox{for some}\,\,\vartheta \in [0, 2
\pi]\}}.$$ Given $\vartheta \in [0, 2 \pi]$, denote by $H(\vartheta
)$ the open half-space such that $\widehat z \in
\partial H(\vartheta )$, and $\xi (\vartheta )$ outward unit normal vector to $H(\vartheta )$ on
$\partial H(\vartheta )$. In particular, $H(0 )= H$, and
\begin{equation}\label{parallel}
\hbox{$\partial H(\tfrac \pi2 )$ is parallel to $\partial H_\nu$.}
\end{equation}
Define the function $m : [0, 2\pi] \to [0, \infty)$ as
%
%
%
%
%
%Let us consider $q\in \partial H_\nu\cap s$, which certainly exists
%since  $\partial H_\nu\cap s\ne\emptyset$ otherwise $\hh
%(\Pi_\nu(\Omega^1))\ge\hh(\Omega^1\cap\partial(H_\nu\cap\Omega^1))$.
%Let now $p$ be a point on $s$ such that $p\ne q$, $p\ne x_1$ and
%$p\ne x_2$. Moreover let $\mu$ be the unit normal, which points
%towards $\Omega^1$ to the support hyperplane containing $s$. If
%$\mathcal{C}=\{\gamma \in\mathbb{S}^{n-1}: \gamma\in
%\text{span}(\mu,\nu)  \}$ we denote by $\gamma(\theta)$ (with
%$\theta\in[0,2\pi]$) the elements of $\mathcal{C}$ and we assume
%that $\gamma(0)=\mu$. With $\pi(\theta)$ we denote the half space
%such that its boundary contains $q$ and is orthogonal. The function
$$m(\vartheta)=\hh(\partial\Omega^1\cap H(\vartheta))\quad \hbox{for $\vartheta  \in [0, 2\pi]$,}$$
and observe that $m$ is a continuous function satisfying
$${m(\pi )=0, \quad m(0)>\frac{P(\Omega ^1)}2.}$$
Therefore there exists $\bar\vartheta \in (0, \pi)$ such that
$m(\bar\vartheta)=\frac{P(\Omega ^1)}2$.
 Since
$\widehat y\ne \widehat z$ we have that $\bar \vartheta\ne \frac\pi 2$,
otherwise
 $H_\nu$ and $H(\tfrac \pi2 )$ would be distinct  half-spaces (since they intersect $\Sigma$ at different
 points),  satisfying
\eqref{parallel} and $\hh (\partial \Omega^1 \cap H_\nu) = \hh (
\partial \Omega ^1 \cap H(\tfrac \pi2 )) = \frac{P(\Omega ^1)}2$. Since $H(\bar \vartheta
)$ fulfills $\hh (\partial \Omega ^1 \cap H(\bar \vartheta
))=\frac{P(\Omega ^1)}2$, and $\widehat z \in \partial H(\bar
\vartheta )$, {one has that $H(\bar \vartheta ) =
H_{\xi (\bar\vartheta)}$}. Hence, there exists a support hyperplane to
$\Omega ^1$ at $\widehat z$ which is orthogonal to $\partial H(\bar
\vartheta )$, otherwise the first equality in \eqref{equality6}
would fail {for $\nu =\xi (\bar\vartheta)$. Such support
hyperplane to $\Omega ^1$, being
orthogonal to $\partial H(\bar \vartheta )$},  would intersect $\Sigma$ only at $\widehat
z$, but this is impossible, since any support hyperplane to $\Omega
^1$ at a point of $\Sigma$ necessarily contains the whole of
$\Sigma$. The strict convexity of $\Omega ^1$ is thus established.
\par\noindent
By the strict convexity of $\Omega ^1$, and the first equality in
\eqref{equality6},
\begin{equation}\label{I1}
\hh(I_\nu(\Omega^1))= \hh (\partial \Omega ^1 \cap H_\nu ) =
P(\Omega ^1)/2 \quad \hbox{for every $\nu \in \mathbb S^{n-1}$,}
\end{equation}
where $I_\nu(\Omega^1)$ denotes the illuminated portion of $\Omega ^1$,
defined as in \eqref{I}. In particular,
\begin{equation}\label{I2}
\hh(I_\nu(\Omega^1))=\hh(I_{-\nu}(\Omega^1))\quad \hbox{for every $\nu \in \mathbb S^{n-1}$}.
\end{equation}
 Property \eqref{I2} implies, via   \cite[Theorem  5.5.11]{Gr},
that $\Omega^1$ is centrally symmetric. Finally, on calling $B$ the
ball with the same perimeter as $\Omega ^1$, we infer from the
second equality in \eqref{equality6}  that
$$\hh(\Pi_\nu(\Omega^1)^+)=\hh(\Pi_\nu(B)^+)\quad \hbox{for every $\nu \in \mathbb S^{n-1}$}.$$
Hence, owing to  \cite[Theorem 5.5.6]{Gr}, we conclude that $\Omega
^1$ is a ball. \qed

%
%
%
%
%
%
%
%
%
%
%
%
%
%
%
%
%$\partial \pi(\bar\theta)$ and $\partial H_\nu$ are parallel
%hyperplanes, both bisecting the perimeter of $\partial\Omega$ and
%they do not coincide (they intersect $s$ in different points). Since
%$\partial\pi(\bar\theta)$ is bisecting and contains $p$, there
%exists a support hyperplane at $p$ which is orthogonal to $\partial
%\pi(\bar\theta)$ otherwise the equality \eqref{equality6} does not
%hold. A contradiction arise since any supporting hyperplanes to $p$
%contains the whole segment $s$.
%
%\textcolor{red}{The strict convexity of $\Omega^1$ and \eqref{equality6} yield $\hh(I_\nu(\Omega^1))=\hh(I_{-\nu}(\Omega^1))$ for all $\nu\in\mathbb{S}^{n-1}$.
%Then Theorems 5.5.11 in \cite{Gr} imply that $\Omega^1$ is centrally symmetric. Not the other hand the second equality in \eqref{equality6} tells us that $\hh{\Pi_nu{\Omega^1}}=\hh{\Pi_nu{\Omega^1}}$ for every $\nu\in\mathbb{S}^{n-1}$ where $B$ is the ball having same parameter as $\Omega^1$. Applying
%theorem 5.5.6 in \cite{Gr} we get that $\Omega ^1$ is a ball. \qed}

\begin{remark}{\rm  The proof of Theorem
\ref{estimate1} considerably simplifies under the additional
assumption that $\Omega$ is convex. In this case, the hyperplane
$H_\nu$ can be defined via \eqref{Hnu} for every $\nu\in\mathbb
S^{n-1}$, and inequality \eqref{medianinequality} follows from
formula \eqref{cauchy} and the chain \eqref{C1}. The
characterization of balls as the only convex sets yielding the
equality in \eqref{medianinequality} can be established as in the
last part of the above proof, on replacing $\Omega ^1$ just by
$\Omega$.}
\end{remark}

%
%We observe that when
%$\Omega$ is convex inequality \eqref{cauchyineq} holds as an equality and indeed it is the well-known \cite[pag.141]{BuZa} Cauchy surface area
%formula
%\begin{equation}\label{cauchy}
%P(\Omega)=\frac1{v_{n-1}}\int_{\mathbb S^{n-1}} {\mathcal
%H}^{n-1}(\Pi_\nu\Omega)\,d\nu.
%\end{equation}
%The open half-space $H_\nu$ %the open half-space, with boundary having
%%outer normal vector $\nu$, such that
%satisfies ${\mathcal H}^{n-1}(\partial \Omega \cap \partial(H_\nu\cap
%\Omega))=P(\Omega)/2$ for all $\nu\in\mathbb{S}^{n-1}$, and once we set
%$$h(\nu)={\mathcal H}^{n-1}(\Omega\cap\partial(H_\nu\cap \Omega))$$
%we have
%\begin{equation}\label{cauchy1}
%P(\Omega)\ge\frac1{v_{n-1}}\int_{\mathbb S^{n-1}} h(\nu)\,d\nu.
%\end{equation}
%Inequalities in \eqref{C1} hold as well and Theorem \ref{estimate1} follows at once.
%%On the other hand, from the definition of $C(\Omega)$ and using
%%\eqref{cauchy1}, we have
%%$$C(\Omega)\ge \frac{P(\Omega)}2\frac1{\displaystyle\min_{H \text{ half-space}\atop
%%{\mathcal H}^{n-1}(\partial\Omega\cap\partial(H\cap \Omega))=P(\Omega)/2} {\mathcal H}^{n-1}(\Omega\cap\partial(H\cap \Omega))}\ge
%%\frac{P(\Omega)}2\frac1{\displaystyle\frac1{\omega_{n}}\int_{\mathbb
%%S^{n-1}} h(\nu)\,d\nu}\ge\frac{\omega_n}{2v_{n-1}}.
%%$$
%\end{remark}}

\smallskip
\par\noindent
We conclude with a proof of Theorem \ref{estimatemean}.

\smallskip
\par\noindent

%
%\begin{theorem}\label{meantrace} Assume that $\Omega$ satisfies
%\eqref{admiss}. Define
%\begin{equation}\label{admiss1}
%C(\Omega)= \frac 2{P(\Omega
%)}\sup_{E\subset\Omega}\frac{P_{C\Omega} (E) \,
%P_{C\Omega}(\Omega\setminus E) }{P_{\Omega} (E)}.
%\end{equation}
%Then $C(\Omega)$ is the best constant in the inequality
%\begin{equation}\label{tracemean}
%\int_{\partial\Omega} |u^* - u^*_{\partial \Omega }| \,d{\mathcal
%H}_{n-1} \le C(\Omega) \|Du \|(\Omega)
%\end{equation}
%for $u \in BV(\Omega )$. Moreover,
%\begin{equation}\label{tracemeanconst}
%C(B) = \begin{cases}  \frac{\omega_n}{2v_{n-1}} & \hbox{if $n \geq
%3$,}
%\\ 2  & \hbox{if $n =2$.}
%\end{cases}
%\end{equation}
%\end{theorem}
%
%%By the same argument as in the proof of Theorem \ref{trace} one
%%can prove the following
%
%
%\begin{theorem}\label{estimatemean}
%Assume that $n \geq 2$. Under the same hypotheses on $\Omega $ as
%in Theorem \ref{mazya},
%\begin{equation}\label{estimatemean1}C(\Omega)\ge C(B) =
%\begin{cases}  \frac{\omega_n}{2v_{n-1}} & \hbox{if $n \geq
%3$,}
%\\ 2  & \hbox{if $n =2$.}
%\end{cases}
%%
%%\frac{\omega_n}{2v_{n-1}} \,\,(= K(B)).
%\end{equation}
%\textcolor{blue}{If $n \geq 3$, then equality holds in
%\eqref{estimatemean1} if and only if $\Omega$ is a ball. Instead,
%if $n=2$, then there exist   open sets $\Omega $, that are not
%disks, for which nevertheless equality holds in
%\eqref{estimatemean1}}.
%\end{theorem}

\par\noindent
{\bf Proof of Theorem \ref{estimatemean}}. Assume that $n \geq 3$.
From \eqref{medmv}
%
%Since $\frac{2s(a-s)}{a} \geq \min\{s, a-s\}$ if $0\leq s \leq a$,
%from \eqref{medianinequality}
 we deduce that $$C_{\rm mv}(\Omega)\ge C_{\rm med}(\Omega )\ge \sqrt{\pi} \,\frac n2 \frac{\Gamma (\frac
{n+1}2)} {\Gamma (\frac {n+2}2)},
%
%\frac{\omega_n}{2v_{n-1}},
$$ namely \eqref{estimatemean1}. Moreover, the
assertion concerning the case of  equality  in \eqref{estimatemean1}
follows from Theorem \ref{estimate1}, since the equality in
\eqref{estimatemean1} implies the equality \eqref{medianinequality}.
%Arguing as in the proof of Theorem \ref{estimate1} we obtain an inequality similar
%to \eqref{C1}, that is,
%\begin{eqnarray*}
%K(\Omega)&\ge& \frac 2{P(\Omega )}
%\bigg(\frac{P(\Omega)}2\bigg)^2\frac1{\displaystyle\inf_{H \text{ half-space}\atop
%{\mathcal H}^{n-1}(\partial\Omega\cap\partial(H\cap \Omega))=P(\Omega)/2}
%{\mathcal H}^{n-1}(\Omega\cap\partial(H\cap \Omega))}\\
%\ge
%\frac{P(\Omega)}2\frac1{\displaystyle\frac1{\omega_{n}}\int_{\mathbb
%S^{n-1}} h(\nu)\,d\nu}\ge\frac{\omega_n}{2v_{n-1}},$$ namely
%\eqref{estimatemean1}.
%As regards the equality case in \eqref{estimatemean1}, the same arguments used in the proof of Theorem \ref{estimate1} imply that $\widetilde\Omega$ (the set of points of density 1 with respect to $\Omega$) is a ball.
\par\noindent
Assume now that $n =2$. Denote by $\widehat \Omega$ the convex hull
of $\overline{\Omega}$. By a standard result in the theory of convex
bodies, there exists an extreme point $x_0$ for $\widehat \Omega$,
and, necessarily, $x_0 \in
\partial \widehat \Omega \cap \partial \Omega$.
%
% $x_0 \in
%\partial \widehat
%\partial \Omega$ for $\Omega$ (in fact, there exist at least $3$ extre
%points of $\widehat \Omega$, which do not lie on the same straight
%line).
Moreover, as a consequence, for instance, of
 Straszewicz's theorem \cite[Theorem
1.4.7]{Schneider}, the point $x_0$ can be chosen in such a way that
it is also
%
%there  exists an extremum point of $\widehat \Omega$ which is also
an exposed point for $\widehat \Omega$, so that
%namely a point $x_0 \in
%\partial \widehat \Omega$  having the property that
there exists an open half-plane $H_0$ such that
$$\Omega \subset H_0,$$
and
$$\overline \Omega
 \cap \partial {H_0} = \{x_0\}.$$
%By another classical result from the theory of convex bodies,
%since $x_0$ is an extremum point of $\widehat \Omega$,
%there exists an extremum point $x_0 \in \partial \widehat
%\Omega$ for $\Omega$ (in fact, there exist at least $3$ extremum
%points of $\widehat \Omega$, which do not lie on the same straight
%line). Another result from the theory of convex bodies ensures
%that
%$x_0 \in
%\partial \widehat \Omega \cap \partial \Omega$.
%By definition of
%extremum point, $x_0$ is not an interior point of any segment
%contained in $\partial \widehat \Omega$.
Assume, without loss of generality, that $x_0= (0,0),$ and $H_0 =
\{(x,y)\in \mathbb R^2: y >0\},$  whence
$$\Omega \subset \{(x,y)\in \mathbb R^2: y >0\}.$$
Given $\varepsilon >0$, consider the open set
$$\Omega (\varepsilon )= \{(x,y)\in \Omega: y < \varepsilon \}.$$
We claim that
\begin{equation}\label{i}
{\mathcal H}^{1}(
 \partial ^M \Omega (\varepsilon) \cap \Omega)
\leq {\mathcal H}^{1}( \partial ^M\Omega (\varepsilon) \cap
\partial  \Omega )
%
%P_\Omega (\Omega (\varepsilon) ) \leq P_{C\Omega} (\Omega
%(\varepsilon) )
\quad \hbox{for $\varepsilon >0$,}
\end{equation}
and
\begin{equation}\label{ii}
\lim _{\varepsilon \to 0} {\mathcal H}^{1}(\partial \Omega \setminus
\partial ^M\Omega (\varepsilon) ) = P(\Omega ).
%
%
%
%P_{C\Omega} (\Omega \setminus \Omega (\varepsilon) ) = P(\Omega ).
\end{equation}
 Let us prove \eqref{i} first. Recall that $\mathcal H^1(\partial \Omega \setminus \partial ^M \Omega) =0$, since $\Omega$ is an admissible domain.
 %One has that
% \begin{equation}\label{10}
% \partial ^M \Omega (\varepsilon) \subset \partial ^M \Omega \cap
% \{y=\varepsilon\}\quad \hbox{for $\varepsilon >0$,}
%\end{equation}
It is easily seen that
\begin{equation}\label{10bis}
 \partial ^M \Omega (\varepsilon) \cap \Omega =
 \{y=\varepsilon\} \cap \Omega \quad \hbox{for $\varepsilon >0$.}
\end{equation}
% Denote the orthogonal projection on $\{y=0\}$ (identified with
%$\R$) by $\Pi$.
\par\noindent
Set $e_2 = (0,1)$. By the coarea formula \eqref{coarea},
\begin{equation}\label{coarea2}
\int _{\partial ^M \Omega (\varepsilon) }|\nu ^\Omega (x) \cdot
e_2|\,d{\mathcal H}^{1}(x) = \int _{\R} {\mathcal H}^{0}((\partial
^M \Omega (\varepsilon) )_z^{e_2})\, dz\,.
\end{equation}
We have that
\begin{align}\label{11}
\int _{\R} {\mathcal H}^{0}((\partial ^M \Omega (\varepsilon)
)_z^{e_2})\, dz & \geq 2 {\mathcal H}^{1}(\Pi_{e_2} (\Omega
(\varepsilon))^+) \geq 2  {\mathcal H}^{1}(\{y=\varepsilon\} \cap
\Omega) = 2{\mathcal H}^{1}(\partial ^M \Omega (\varepsilon) \cap
\Omega).
\end{align}
On the other hand,
\begin{align}\label{12}
\int _{\partial ^M\Omega (\varepsilon) }|\nu ^\Omega (x) \cdot
e_2|\,d{\mathcal H}^{1}(x)
 & \leq {\mathcal H}^{1}(\partial ^M\Omega (\varepsilon) ) \\ \nonumber &
= {\mathcal H}^{1}( \partial ^M\Omega (\varepsilon) \cap \Omega  )
  + {\mathcal H}^{1}( \partial ^M\Omega (\varepsilon)
\cap \partial  \Omega ).
%\\ \nonumber & =
% {\mathcal H}^{1}( \{y=\varepsilon\} \cap \Omega  ) + {\mathcal H}^{1}( \partial ^M\Omega (\varepsilon) \cap \partial ^M \Omega )
%\\ \nonumber & =
%{\mathcal H}^{1}(
% \partial ^M \Omega (\varepsilon) \cap \Omega)
%+ {\mathcal H}^{1}( \partial ^M\Omega (\varepsilon) \cap \partial
%^M \Omega )\,.
\end{align}
Inequality \eqref{i} follows from \eqref{11} and \eqref{12}.
\par\noindent
Consider next \eqref{ii}. Owing to our choice of $x_0$, one has that
$\partial ^M\Omega \setminus \partial ^M\Omega (\varepsilon) $ is an
increasing family of sets as $\varepsilon \searrow 0^+$. Moreover,
\begin{equation}\label{14}
 \cup_{\varepsilon >0} \, (\partial ^M\Omega \setminus \partial ^M\Omega (\varepsilon) )=  \partial ^M \Omega \setminus \Gamma\,,
\end{equation}
where $\Gamma$ equals either $\{x_0\}$ or $\emptyset$, according to
whether $x_0$ belongs to $\partial ^M \Omega$ or not. In order to
verify \eqref{14}, observe that
$$\partial ^M \Omega \cap \{y <\varepsilon\} \subset \partial ^M
\Omega \cap \partial ^M (\Omega \cap \{y <\varepsilon\})= \partial
^M \Omega \cap \partial ^M \Omega (\varepsilon ),$$ where the
inclusion holds by \eqref{dic10}, with $H_{\nu ,i}^{\pm}$ replaced
with $\{y <\varepsilon\}$. Thus,
$$\partial ^M
\Omega \cap \partial ^M \Omega (\varepsilon ) = \partial ^M \Omega
\setminus ( \partial ^M\Omega  \cap \partial ^M  \Omega (\varepsilon
)) \subset \partial ^M \Omega \setminus (\partial ^M \Omega \cap \{y
<\varepsilon\}).$$ Now, $\partial ^M \Omega \cap \{y <\varepsilon\}
\searrow \Gamma$ as $\varepsilon \to 0^+$, since $\partial ^M \Omega
\subset \partial \Omega$ and $\partial  \Omega \cap \{y
<\varepsilon\} \searrow \{x_0\}$. Hence, \eqref{14} follows.
\par\noindent
 Since ${\mathcal H}^{1}$ is a measure on Borel
sets, from \eqref{14} we obtain that
\begin{equation}\label{D}
%\lim _{\varepsilon \to 0^+} P_{C\Omega}(\Omega \setminus \Omega (\varepsilon) ) =
\lim _{\varepsilon \to 0^+} {\mathcal H}^{1}(\partial ^M \Omega \setminus
\partial ^M\Omega (\varepsilon) ) =  {\mathcal H}^{1}(\partial ^M \Omega \setminus \Gamma ) =  {\mathcal H}^{1}(\partial ^M \Omega ).
\end{equation}
{Inasmuch as $\Omega$ is an admissible domain,
equation \eqref{ii} is a consequence of \eqref{D}.}
%
%
%
%Hence,
%$$\lim _{\varepsilon \to 0} {\mathcal H}^{1}(\partial^M \Omega \setminus
%\partial ^M\Omega (\varepsilon) ) = {\mathcal H}^{1}(\partial ^M
%\Omega),$$
%and \eqref{ii}
%is thus established.
\par\noindent
By \eqref{i} and \eqref{ii}
\begin{eqnarray}\label{15}
C_{\rm mv}(\Omega)   &\geq& \lim _{\varepsilon \to 0^+} \frac 2{{\mathcal
H}^{1}(\partial  \Omega )}\frac{{\mathcal H}^{1}( \partial
\Omega (\varepsilon) \cap \partial ^M \Omega ) \, {\mathcal
H}^{1}(\partial ^M \Omega \setminus
\partial ^M\Omega (\varepsilon) )}{{\mathcal H}^{1}(
 \partial ^M \Omega (\varepsilon) \cap \Omega)} \\ \nonumber &\geq& \lim _{\varepsilon \to 0^+} \frac {2\,{\mathcal
H}^{1}(\partial \Omega \setminus
\partial ^M\Omega (\varepsilon) ) }{{\mathcal H}^{1}(\partial
\Omega )}  = 2\,,
\end{eqnarray}
whence \eqref{estimatemean2} follows.
\par
We complete the proof by showing that there exist domains $\Omega$
in  $\mathbb R^2$, with a continuously differentiable boundary,
which are different from a disk, such that
%\begin{equation}\label{stadium}
$$C_{\rm mv}(\Omega)=2.$$
%\end{equation}
The sets $\Omega$ which will be exhibited are stadium-shaped.
{Let us preliminarily observe that,  if $\Omega$ is a
bounded convex domain in $\mathbb R^2$, then
\begin{equation}\label{isopconstantpc}
C_{\rm mv}(\Omega) = \frac 2{\hu (\partial \Omega )} \sup
_{\overset{E={H}\cap\Omega}{{H} \>{\text{half-plane}}
%\atop{\partial E}\cap\Omega\> \text{has only one connected component}
} } \frac{\hu (\partial E \cap
\partial \Omega ) \,\, \hu ( \partial \Omega \setminus \partial E)}{\hu
(\partial E \cap \Omega)}.
\end{equation}}
Indeed, as a consequence of the results of \cite[Section
9.4.1]{Mazbook},
\begin{equation}\label{isopconstantp}
C_{\rm mv}(\Omega)= \frac 2{\hu (\partial \Omega )} \sup _{
 \overset
{E={\cal P}\cap\Omega}{
%}\atop{
\cal P \>{\text{polygon}}}} \frac{\hu (\partial E \cap
\partial \Omega ) \,\, \hu ( \partial \Omega \setminus \partial E)}{\hu
(\partial E \cap \Omega)}.
\end{equation}
%where the supremum is taken over polygons in the class of admissible
%sets.
% We claim that  in \eqref{isopconstantp} the supremum can be
%restricted to the class of those sets $E \subset \Omega$ which are
%the intersection of an half-plane with $\Omega$, namely
%%\begin{equation}\label{isopconstantpc}
%%K(\Omega ) = \frac 2{\hh (\partial \Omega )} \sup _{{E={\cal P}\cap\Omega}\atop{{\cal P \>{\text{polygon}}}\atop{\partial E}\cap\Omega\> \text{has only one connected component}}}} \frac{\hh
%%(\partial ^M E \cap
%%\partial \Omega ) \,\, \hh ( \partial \Omega \setminus \partial ^M E)}{\hh
%%(\partial ^M E \cap \Omega)}.
%%\end{equation}
%\begin{equation}\label{isopconstantpc}
%C_{\rm mv}(\Omega) = \frac 2{\hu (\partial \Omega )} \sup _{\overset{E={H}\cap\Omega}{{H} \>{\text{half-plane}}
%%\atop{\partial E}\cap\Omega\> \text{has only one connected component}
%}
%}
%\frac{\hu
%(\partial E \cap
%\partial \Omega ) \,\, \hu ( \partial \Omega \setminus \partial E)}{\hu
%(\partial E \cap \Omega)}.
%\end{equation}
To prove this claim, observe that, if $\{E_i\}_{i=1,\dots,k}$,
$k\in\N$, are the connected components of a set $E$ as in
\eqref{isopconstantp}, then
\begin{align}\label{E}
\frac{\hu (\partial E \cap
\partial \Omega ) \,\, \hu ( \partial \Omega \setminus \partial E)}{\hu
(\partial E \cap \Omega)}  & =
 \frac{\hu(\partial (\cup _{i=1}^k E_i) \cap
\partial \Omega ) \,\, (\hu ( \partial \Omega \setminus \partial (\cup _{i=1}^k E_i)))}
{\hu (\partial (\cup _{i=1}^k E_i) \cap \Omega)}
\\ \nonumber
&  = \frac{\hu( \cup _{i=1}^k (\partial E_i) \cap
\partial \Omega )
\,\, (\hu ( \partial \Omega \setminus (\cup _{i=1}^k \partial
E_i)))}{\hu (\cup _{i=1}^k (\partial E_i) \cap \Omega)}
\\ \nonumber
&  =\frac{(\sum_{i=1}^k\hu (\partial E_i \cap
\partial \Omega )) \,\, (\hu ( \partial \Omega )-\sum_{i=1}^k\hu
(\partial E_i \cap
\partial \Omega ))}{\sum_{i=1}^k\hu
(\partial E_i \cap \Omega)}
\\ \nonumber
& \le\max_i  \frac{\hu (\partial E_i \cap
\partial \Omega ) \,\,( \hu ( \partial \Omega )-\hu
(\partial E_i \cap
\partial \Omega ))}{\hu
(\partial E_i \cap \Omega)}
\\ \nonumber
& =\max_i  \frac{\hu (\partial E_i \cap
\partial \Omega ) \,\, \hu ( \partial \Omega \setminus \partial E_i )}{\hu
(\partial E_i \cap \Omega)} .
\end{align}
Thus, on replacing, if necessary, $E$ with one of its connected components, the supremum in \eqref{isopconstantp} can be restricted to
the class of sets $E$ of the form ${\cal P}\cap\Omega$ which are
connected. Our next step consists in showing  that we may also assume
that  $\partial E \cap \Omega$ is connected. Indeed, given any
connected set of the form  $E={\cal P}\cap\Omega$ for some polygon
${\cal P}$, denote by $\{F_j\}_{j=1,\dots,m}$, $m\in\N$, the
connected components of $\Omega \setminus E$.
Since $F_j$ and
$\Omega \setminus F_j$ are connected for every $j=1,\dots,m$, we
have that $\partial F_j \cap \Omega$ is connected as well.
 Equation \eqref{boundaries4}, and an analogous chain as in \eqref{E} tells us that
\begin{eqnarray*}
\frac{\hu
(\partial E \cap
\partial \Omega ) \,\, \hu ( \partial \Omega \setminus \partial E)}{\hu
(\partial E \cap \Omega) } &= &  \frac{\hu (\partial (\Omega
\setminus E) \cap
\partial \Omega ) \,\, \hu ( \partial \Omega \setminus \partial (\Omega \setminus E))}{\hu
(\partial (\Omega \setminus E) \cap \Omega)}
\\
& = & \frac{\hu (\partial (\Omega \setminus E) \cap
\partial \Omega ) \,\,( \hu ( \partial \Omega )-\hu
(\partial (\Omega \setminus E)  \cap
\partial \Omega ))}{\hu
(\partial (\Omega \setminus E) \cap \Omega)}\\ &\le &\max_j
\frac{\hu (\partial F_j \cap
\partial \Omega ) \,\,( \hu ( \partial \Omega )-\hu
(\partial F_j \cap
\partial \Omega ))}{\hu
(\partial F_j \cap \Omega)
}
\\
&= & \max_j  \frac{\hu (\partial F_j \cap
\partial \Omega ) \,\, \hu ( \partial \Omega \setminus \partial F_j )}{\hu
(\partial F_j \cap \Omega)
}
.
\end{eqnarray*}
This proves that, on further replacing, if necessary, any connected
set $E$ of the form ${\cal P}\cap\Omega$ with one of the connected
components of $\Omega \setminus E$, the supremum in
\eqref{isopconstantp} can be restricted to the class of sets $E$ of
the form ${\cal P}\cap\Omega$, which are connected and such that
$\partial E \cap \Omega$ is also connected. Since $\Omega$ is
convex, the relevant supremum can be finally restricted to the class
of sets $E$  such that $\partial E \cap \Omega$ is a straight
segment. Hence, our claim follows.
\par
%
% over sets $E$
%such that $E$ and $\Omega\setminus E$ are connected. Eventually
%$\partial E \cap \Omega$ is a segment.
\par\noindent
Now, define
%in a position to compute the constant $C_{\rm
%mv}(S_{R,d})$, where %$S_{R,d}$ is the stadium
\begin{equation}
S_{R,d}= \text{convex hull of  two disks of equal radii $R$,  with
centers at distance } d,
\end{equation}
a stadium-shaped domain, with semi-perimeter
 $p=d+\pi R$. Let us introduce the curvilinear abscissa $s\in\R$ on
$\partial S_{R,d}$, and the periodic parametrization $\R \ni s
\mapsto x(s)\in \mathbb R^2$ of $\partial S_{R,d}$ with respect to
$s$. Observe that, for any half-plane $H$, the set $\partial H\cap
 S_{R,d}$ (if not empty) is a chord whose terminals split $\partial S_{R,d}$ into two parts of lengths $a$ and $2p-a$, for some $a\in (0, p]$.
Thus, such a chord has terminals $x(s)$ and $ x(s+a)$, for some
$s$, and its length is given by
\begin{equation}\label{elle}
\ell_a(s)=|x(s+a)- x(s)|.
\end{equation}
Hence, owing to \eqref{isopconstantpc},
\begin{equation}\label{ellesup}
C_{\rm mv}(S_{R,d})=\frac1p\sup_{0<a\le
p}\frac{a(2p-a)}{\displaystyle \min_s\ell_a(s)}.
\end{equation}
For each $a\in (0, p]$, the function $\ell_a$ is continuously
differentiable in $\mathbb R$, and periodic of period $p$. Moreover,
 $s$ is a  stationary point of $\ell_a(s)$ if and only if
\begin{equation}\label{ellecond}
( x(s+a)- x(s))\cdot( x'(s+a)- x'(s))=0,
\end{equation}
where $\cdot$ stands for scalar product in $\mathbb R^2$. Condition
\eqref{ellecond} entails that the tangent straight-lines to
$\partial S_{R,d}$ at $ x(s)$ and $ x(s+a)$ are either orthogonal to
the chord $x(s+a)- x(s)$ or they meet at a point which is
equidistant from $x(s)$ and $ x(s+a)$.
%\textcolor{red}{on
%the same side}, form equal angles with the chord.
It is easily verified that \eqref{ellecond} %then implies that any chord which minimizes \eqref{elle}
is  satisfied only if one of the following situations occurs:
\begin{description}
\item[\rm$\>$ (i)] the chord  is parallel to the flat parts of $\partial S_{R,d}$;
\item[\rm$\>$(ii)] the chord is orthogonal to the flat parts of $\partial S_{R,d}$;
\item[\rm(iii)] the chord has terminals  belonging to the same half circle of $\partial S_{R,d}$.
\end{description}
Elementary arguments show that, in case (i), the chord cannot be a
minimizer of $\ell_a(s)$, whereas a minimizer with property (iii)
can always be chosen in such a way that it is orthogonal to the flat
parts of $\partial S_{R,d}$.
\par\noindent In conclusion, in order to seek for the minimum of
$\ell_a(s)$ in $s$, we may restrict our analysis  to those values of
$s$ such that the chord with terminals  $x(s)$ and $ x(s+a)$ is
 orthogonal to
the flat parts of $\partial S_{R,d}$.  It is then easily seen that
\begin{equation*}
%g(a)=
\min_s\ell_a(s)= \left\{
\begin{array}{ll}
2R& \hbox{if \,\,} a\ge \pi R\\
2R\sin\frac a{2R}& \hbox{if \,\,}a< \pi R.
\end{array}\right.
\end{equation*}
%and
%\begin{equation}\label{ellesup}
%K(S_{R,d})=\frac1p\sup_{0<a\le p}\frac{a(2p-a)}{g(a)}.
%\end{equation}
A straightforward computation now shows that, if $0<d\le (4-\pi)R$,
then the supremum in \eqref{ellesup} is achieved in the limit as $a$
goes to 0.  Hence, $C_{\rm mv}(S_{R,d})=2$.
%, if $0<d\le \frac{4-\pi}2R$, then $K(S_{R,d})=2$
%\begin{equation}\label{kstadio}
%K(S_{R,d})=\frac1p\lim_{a\rightarrow0}\frac{a(2p-a)}{g(a)}=2, \qquad \text{if }0<d\le (4-\pi)R.
%\end{equation}
\qed

\begin{remark}
{\rm An inspection of the  proof of Theorem \ref{estimatemean}
shows that
there do exist sets $\Omega$ for which $C_{\rm mv}(\Omega)>2$. Indeed, one can verify that %instead of \eqref{kstadio},
\begin{equation*}%\label{kstadio}
C_{\rm mv}(S_{R,d})=\frac{d+\pi R}{2R}>2,
\end{equation*}
provided that $d> (4-\pi)R$. }
\end{remark}

 {}{}{}{}{}{}
 {}{}{}{}{}{}{}{}{}{}{}{}{}{}{}{}{}{}{}{}{}{}{}{}{}{}{}{}{}{}{}{}{}{}{}{}{}{}
 {}{}{}{}{}{}
 {}{}{}{}{}{}{}{}{}{}{}{}{}{}{}{}{}{}{}{}{}{}{}{}{}{}{}{}{}{}{}{}{}{}{}{}{}{}

\bigskip
\bigskip
\par\noindent
{\bf Acknowledgement}. The authors wish  to thank Bernhard Kawohl
for several stimulating discussions. \par\noindent This research was
partly supported by the research project  Prin 2008 ``Geometric
aspects of partial differential equations and related topics" of
MIUR (Italian Ministry of University).

\end{document}